\documentclass[11pt]{article}

\usepackage{listings}
\usepackage{pdflscape}
\usepackage{epsfig}
\usepackage{lscape}
\usepackage{longtable}
\usepackage{amsmath,amssymb}
\usepackage{graphicx}
\usepackage{color}
\usepackage{placeins}
\usepackage{url}
\def\sqr#1#2{{\vcenter{\vbox{\hrule height.#2pt
        \hbox{\vrule width.#2pt height#1pt \kern#2pt
        \vrule width.#2pt}
        \hrule height.#2pt}}}}
\newcommand{\qed}{\hfill $ \mathchoice\sqr74\sqr74\sqr{2.1}3\sqr{1.5}3 $}
\def\approxleq{ \kern3pt \mbox{\raisebox{.6ex}{$<$}} \kern-8pt
  \mbox{\raisebox{-.6ex}{$\sim$}} \kern5pt}

\def\inprod#1#2{\langle#1,\,#2 \rangle}
\def\Inprod#1#2{\left\langle#1,\,#2 \right\rangle}

 \def\cU{{\cal U}}    
   \def\cW{{\cal W}}
     \def\cV{{\cal V}}
\def\cI{{\cal I}}

\def\disp{\displaystyle}

\newlength{\len}
\newtheorem{theorem}{Theorem}[section]
\newtheorem{prop}{Proposition}[section]

\newtheorem{remark}{Remark}[section]
\newtheorem{assumption}{Assumption}[section]

\begin{document}
\title{On the convergence properties of a majorized ADMM for linearly constrained convex optimization problems with coupled objective functions}

\author{Ying Cui\thanks {Department of Mathematics, National University of Singapore, 10 Lower Kent Ridge Road, Singapore ({\tt cuiying@nus.edu.sg}).}, \;Xudong Li\thanks{Department of Mathematics, National University of Singapore, 10 Lower Kent Ridge Road, Singapore ({\tt lixudong@nus.edu.sg}).}, \; Defeng Sun\thanks{Department  of  Mathematics  and  Risk  Management  Institute, National University of Singapore, 10 Lower Kent Ridge Road, Singapore ({\tt matsundf@nus.edu.sg}). } \  and \ Kim-Chuan Toh\thanks{Department of Mathematics, National University of Singapore, 10 Lower Kent Ridge Road, Singapore
({\tt mattohkc@nus.edu.sg}).
 }\\[10pt]
\textit {Dedicated to Lucien Polak on the occasion of his 85th birthday}
}

\maketitle

\begin{abstract}
In this paper, we establish the  convergence properties for a  majorized alternating direction method of multipliers (ADMM) for linearly constrained convex optimization problems whose objectives contain coupled functions.
Our convergence analysis relies on the
generalized Mean-Value Theorem which plays an important role  to properly control the cross terms due to the presence of coupled objective functions.
Our results in particular show that directly applying 2-block ADMM {with a large step length} to the  linearly constrained convex optimization problem with a quadratically coupled objective function is convergent under mild conditions.
We also provide several iteration complexity results for the  algorithm.

\end{abstract}

\medskip
{\small
\begin{center}
\parbox{0.95\hsize}{{\bf Keywords.}\; ADMM, Coupled objective function, Convex quadratic programming, Majorization, Iteration complexity, Nonsmooth analysis
}
\end{center}
}

\section{Introduction}\label{Introduction}
Consider the following convex optimization problem:
\begin{equation}\label{opt}\begin{array}{ll}
\displaystyle\min_{u,v}& \theta(u,v): =  p(u) + q(v) + \phi(u,v),\\ [8pt]
\text{s.t.} & \mathcal{A}^*u + \mathcal{B}^*v = c,
\end{array}
\end{equation}
where $p: \mathcal{U}\to(-\infty,\infty]$, $q: \mathcal{V}\to(-\infty,\infty]$ are two closed proper convex functions (possibly nonsmooth), $\phi:\mathcal{U}\times \mathcal{V}\to (-\infty, \infty)$ is a smooth convex function whose gradient mapping is  Lipschitz continuous,
$\mathcal{A}: \mathcal{X}\to\mathcal{U}$ and
$\mathcal{B}:\mathcal{X}\to\mathcal{V}$ are two given linear operators, $c\in\mathcal{X}$ is a given vector,  and $\mathcal{U},\mathcal{V}$ and  $\mathcal{X}$ are three real finite dimensional Euclidean spaces each equipped with an inner product $\langle \cdot,\cdot\rangle$ and its induced norm $\|\cdot\|$.

Many interesting optimization problems are of the form (\ref{opt}).
One particular case is the following problem whose objective is the sum of a quadratic function and a squared distance function to a closed convex set:
\begin{equation}\label{ex}\begin{array}{ll}
\min & \displaystyle\frac{1}{2}\bigg\langle \left(\begin{array}{c} u\\v\end{array}\right), \widetilde{\mathcal{Q}}\left(\begin{array}{c} u\\v\end{array}\right)\bigg\rangle + \frac{\rho}{2}\bigg\|\left(\begin{array}{c} u\\v\end{array}\right) - \Pi_{{\mathcal{K}_1}}\left(\begin{array}{c} u\\v\end{array}\right)\bigg\|^2, \\[8pt]
\text{s.t.} & \mathcal{A}^*u + \mathcal{B}^*v = c, \\[8pt]
&
u\in \mathcal{K}_2, ~v\in{\mathcal{K}_3},
\end{array}
\end{equation}
where $\rho>0$ is a penalty parameter, $\widetilde{\mathcal{Q}}:\mathcal{U}\times\mathcal{V}\to\mathcal{U}\times \mathcal{V}$ is a self-ajoint positive semidefinite linear operator,  $\mathcal{K}_1\subseteq \mathcal{U}\times\mathcal{V}$, $\mathcal{K}_2\subseteq \mathcal{U}$ and $\mathcal{K}_3\subseteq \mathcal{V}$ are closed convex sets and $\Pi_{\mathcal{K}_1}(\cdot, \cdot)$ denotes the metric projection onto $\mathcal{K}_1$.

One popular way to solve problem (\ref{opt}) is the augmented Lagrangian method (ALM).
Given the Lagrangian multiplier $x\in\mathcal{X}$ of the linear constraint in (\ref{opt}), the augmented Lagrangian function associated with the parameter $\sigma >0$ is defined as
\begin{equation}
\mathcal{L}_\sigma(u,v;x) = \theta(u,v) + \langle x, \mathcal{A}^*u + \mathcal{B}^*v - c\rangle +\frac{\sigma}{2}\|\mathcal{A}^*u + \mathcal{B}^*v -c\|^2, \quad (u,v)\in\mathcal{U}\times\mathcal{V}.
\end{equation}
The ALM minimizes $\mathcal{L}_\sigma(u,v;x)$ with respect to $(u,v)$ simultaneously regardless of whether the objective function is coupled or not before updating the Lagrangian multiplier $x$ along the gradient ascent direction. Numerically, however, to minimize $\mathcal{L}_\sigma(u,v;x)$ with respect to $(u,v)$ jointly may be a difficult task due to the non-separable structure of $\theta(\cdot,\cdot)$ combined with the nonsmoothness of $p(\cdot)$ and $q(\cdot)$.

When  the objective function in (\ref{opt}) is separable for $u$ and $v$, one can alleviate the numerical difficulty in the ALM by directly applying the  alternating direction method of multipliers (ADMM). The iteration scheme of the ADMM works as follows:
\begin{equation}\label{alm}
\left\{\begin{array}{ll}
u^{k+1} =\displaystyle\arg\min_{u} \mathcal{L}_\sigma(u,v^k; x^k),\\[8pt]
v^{k+1} = \displaystyle\arg\min_{v} \mathcal{L}_\sigma(u^{k+1},v; x^k), \\[8pt]
x^{k+1}  = x^k + \tau\sigma(\mathcal{A}^*u^{k+1} + \mathcal{B}^*v^{k+1} - c),
\end{array}\right.
\end{equation}
where $\tau>0$ is the step length. The global convergence of the  ADMM with $\tau\in(0,\frac{1+\sqrt{5}}{2})$ and a separable objective function has been extensively studied in the literature, see, for examples, ~\cite{Gabay83, Gabay1976, Glowinski1980lectures, Glowinski1975, Eckstein_1992}. For a recent survey, see Eckstein and Yao~\cite{Eckstein2014}. Although it is possible to apply the ADMM directly to problem (\ref{opt}) even if $\phi(\cdot,\cdot)$ is not separable, its convergence analysis is largely non-existent. One way to deal with the non-separablity of $\phi(\cdot,\cdot)$ is to introduce a new variable $w = \left(\begin{array}{c}u\\ v\end{array}\right)$. By letting $\widetilde{\mathcal{A}} = \left(\begin{array}{c}\mathcal{A}\\ {\cI}_1\\0\end{array}\right)$, $\widetilde{\mathcal{B}} = \left(\begin{array}{c}\mathcal{B}\\ 0\\{\cI}_2\end{array}\right)$, $\widetilde{\mathcal{C}} = \left(\begin{array}{cc}0 & 0\\ {\cI}_1 & 0\\0 &{\cI}_2\end{array}\right)$ and $\tilde{c} = \left(\begin{array}{c} c\\0 \\0 \end{array}\right)$ {with identity maps $\cI_1:\cU \to \cU$ and $\cI_2:\cV\to\cV$}, we can rewrite the optimization problem (\ref{ex}) equivalently as
\begin{equation}\label{opt_3}\begin{array}{ll}
\displaystyle\min_{u,v,w} & \tilde{\theta}(u,v,w): = p(u)  + q(v) + \phi(w),\\[8pt]
\text{s.t.} & \widetilde{\mathcal{A}}^*u + \widetilde{\mathcal{B}}^*v  + \widetilde{\mathcal{C}}^*w = \tilde{c}.
\end{array}
\end{equation}
For given $\sigma>0$,  the corresponding augmented Lagrangian function for problem (\ref{opt_3}) is
\begin{equation*}
\widetilde{L}_\sigma(u,v,w; x) = \tilde{\theta}(u,v,w) + \langle x, \widetilde{\mathcal{A}}^*u + \widetilde{\mathcal{B}}^*v + \widetilde{\mathcal{C}}^*w - \tilde{c}\rangle +  \displaystyle\frac{\sigma}{2}\| \widetilde{\mathcal{A}}^*u + \widetilde{\mathcal{B}}^*v + \widetilde{\mathcal{C}}^*w - \tilde{c}\|^2,
\end{equation*}
where $(u,v,w)\in\mathcal{U}\times\mathcal{V}\times(\mathcal{U}\times\mathcal{V})$ and $x\in\mathcal{X}$.
Directly applying the 3-Block ADMM yields the following framework:
\begin{equation*}
\left\{\begin{array}{ll}
u^{k+1} =\displaystyle\arg\min_{u} \widetilde{\mathcal{L}}_\sigma(u,v^k,w^k; x^k),\\[8pt]
v^{k+1} =\displaystyle\arg\min_{v} \widetilde{\mathcal{L}}_\sigma(u^{k+1},v,w^k; x^k),\\[8pt]
w^{k+1} =\displaystyle\arg\min_{w} \widetilde{\mathcal{L}}_\sigma(u^{k+1},v^{k+1},w; x^k),\\[8pt]
x^{k+1}  = x^k  + \tau\sigma(\widetilde{\mathcal{A}}^*u^{k+1} + \widetilde{\mathcal{B}}^*v^{k+1} + \widetilde{\mathcal{C}}^*w^{k+1} - \tilde{c}),\end{array}\right.
\end{equation*}
where $\tau >0$ is the step length.
Even though numerically the 3-block ADMM works well  for many applications, generally it is not a convergent algorithm even if $\tau$ is as small as $10^{-8}$ as shown in the 
 counterexamples given by Chen et al.~\cite{Chen_divergence2013}.

In this paper, we will conduct a thorough convergence analysis about the 2-block ADMM when it is applied to problem (\ref{opt}) with non-separable objective functions. Unlike the case with separable objective functions, there are very few papers on the ADMM targeting the problem (\ref{opt}) except for the work of Hong et al.~\cite{Hong2014}, where the authors studied a majorized multi-block ADMM  for  linearly constrained optimization problems with non-separable objectives. When specialized  to the 2-block case for problem (\ref{opt}), their algorithm works as follows:
\begin{equation}\label{Hong_iter}\left\{\begin{array}{ll}
u^{k+1} =\displaystyle\arg\min_{u}\{ p(u)  +\langle x^k,\mathcal{A}^*u \rangle +
\hat{h}_1(u; u^k, v^k)\},\\[8pt]
v^{k+1} = \displaystyle\arg\min_{v} \{q(v) + \langle x^k, \mathcal{B}^*v\rangle +\hat{h}_2(v; u^{k+1}, v^k)\}, \\[8pt]
x^{k+1}  = x^k + \alpha_k\sigma(\mathcal{A}^*u^{k+1} + \mathcal{B}^*v^{k+1} - c),
\end{array}\right.\end{equation}
where $\hat{h}_{1}(u; u^k,v^k)$ and $\hat{h}_{2}(v; u^{k+1},v^k)$ are  majorization functions of $\phi(u,v) + \frac{\sigma}{2}\|\mathcal{A}^*u + \mathcal{B}^*v - c\|^2$ at $(u^k,v^k)$
and $(u^{k+1}, v^{k})$, respectively and $\alpha_k >0$ is the step length.
Hong et al.~\cite{Hong2014} provided a very general convergence analysis of their majorized ADMM assuming that the step length $\alpha_k$ 
is a sufficiently small fixed number or converging to zero,
among other conditions. Since a large step length is almost always desired in practice, one needs to develop a new convergence theorem beyond the one in~\cite{Hong2014}.
{Similar to} Hong et al.'s work \cite{Hong2014},
our approach also relies on the majorization technique applied to the smooth coupled function $\phi(\cdot,\cdot)$. One difference  is that we majorize $\phi(u,v)$ at $(u^k,v^k)$ before the $(k+1)$th iteration instead of changing the majorization function based on $(u^{k+1},v^k)$ when updating $v^{k+1}$ as in (\ref{Hong_iter}). Interestingly, if $\phi(\cdot,\cdot)$ merely consists of quadratically coupled functions and separable smooth functions,
our majorized ADMM is exactly the same as the one proposed by Hong et al. under a proper choice of the majorization functions.  Moreover,  for applications like (\ref{ex}), a potential advantage of our method is that we only need to compute
the projection $\Pi_{\mathcal{K}_1}(\cdot,\cdot)$ once in order to compute
 $\nabla \phi(\cdot, \cdot)$ as a part of the majorization function  within one iteration, while {the procedure} (\ref{Hong_iter}) needs
 to {compute $\Pi_{\mathcal{K}_1}(\cdot,\cdot)$ at}
 two different points $(u^k,v^k)$ and $(u^{k+1}, v^k)$. In the subsequent discussions one can see that
  by making use of nonsmooth analysis, especially the generalized Mean-Value Theorem, we are  able to establish the global convergence  and the iteration complexity for {our majorized ADMM} with the step length $\tau\in(0,\frac{1+\sqrt{5}}{2})$. To the best of our knowledge, this is the first paper providing the convergence {properties} of {the} majorized ADMM {with a large step length} for solving {linearly constrained convex optimization problems} with coupled smooth objective functions.

%

The remaining parts of our paper are organized as follows. In the next section, we  provide some preliminary results. Section 3 focuses on our framework of a majorized ADMM and two important inequalities for the convergence analysis.
In Section 4, we prove the global convergence and several  iteration complexity results of the proposed algorithm.  We conclude our paper in the last section.

\section{Preliminaries}
In this section, we shall provide some preliminary results  that will be used in our subsequent discussions.

Denote $w:= \left(\begin{array}{c} u\\ v \end{array}\right)$. Since $\phi(\cdot)$ is assumed to be a convex function with a Lipschitz continuous gradient,
 $\nabla \phi(\cdot)$ is globally Lipschitz continuous and $\nabla^2 \phi(\cdot)$ exists almost everywhere.
Thus, the following Clarke's generalized Hessian at given $w\in\mathcal{U}\times\mathcal{V}$ is well defined~\cite{Clarke}:
\begin{equation}\label{gen_Hessian}
\partial^2 \phi(w) = \text{conv}\{\displaystyle\lim_{w^k\to w} \nabla^2 \phi(w^k), \nabla^2 \phi(w^k) \,\text{exists}\},
\end{equation}
where ``conv$\{S\}$'' denotes the convex hall of a given set $S$.
Note that  $\mathcal{W}$ is self-adjoint and positive semidefinite, i.e.,
$\mathcal{W}\succeq 0$, for any $\mathcal{W}\in\partial^2\phi(w)$, $w\in\mathcal{U}\times\mathcal{V}$.
In~\cite{JBHU_1984},  Hiriart-Urruty and Nguyen provide a second order Mean-Value Theorem for $\phi$, which states that
for any $w'$ and $w$ in $\mathcal{U}\times\mathcal{V}$, {there exists  $z\in[w',w]$ and $\mathcal{W}\in\partial^2\phi(z)$ such that}
\begin{equation*}
\phi(w) = \phi(w') + \langle \nabla\phi(w') , w-w'\rangle +\frac{1}{2}\langle w-w', \mathcal{W}(w-w')\rangle,
\end{equation*}
where $[w',w]$ denotes the line segment connecting $w'$ and $w$.

Since $\nabla\phi$ is globally Lipschitz continuous, there exist two self-adjoint positive semidefinite linear operators $\mathcal{Q}$ and $\mathcal{H}: \mathcal{U}\times\mathcal{V}\to\mathcal{U}\times\mathcal{V}$
such that for any $w\in\mathcal{U}\times\mathcal{V}$,
\begin{equation}\label{bound}
\mathcal{Q} \preceq \mathcal{W} \preceq \mathcal{Q}+\mathcal{H} \quad\forall \,\mathcal{W}\in\partial^2 \phi(w).
\end{equation}
Thus, for any $w, w'\in\mathcal{U}\times\mathcal{V}$, we have
\begin{equation}\label{phi_convex}
\phi(w)\geq \phi(w') + \langle \nabla \phi (w'), w - w' \rangle+\frac{1}{2}\|w'-w\|^2_{\mathcal{Q}}
\end{equation}
 and
\begin{equation}\label{phi_majorize}
\phi(w) \leq \hat{\phi}(w; w'): = \phi(w') + \langle \nabla \phi (w'), w - w' \rangle+\frac{1}{2}\|w'-w\|^2_{\mathcal{Q}+\mathcal{H}}.\end{equation}
In this paper we further assume that
\begin{equation}\label{diag_H}
\mathcal{H} = \text{Diag}\,(\mathcal{D}_1, \mathcal{D}_2),
\end{equation}
 where $\mathcal{D}_1: \mathcal{U}\to\mathcal{U}$ and $\mathcal{D}_2: \mathcal{V}\to\mathcal{V}$ are two self-adjoint positive semidefinite linear operators. In fact, this kind of structure naturally appears in applications like (\ref{ex}), where
the best possible lower bound of the generalized Hessian is $\widetilde{\mathcal{Q}}$ and the best possible upper bound of the generalized Hessian is $\widetilde{\mathcal{Q}}+\mathcal{I}$, where $\mathcal{I}: \mathcal{U}\times\mathcal{V}\to\mathcal{U}\times\mathcal{V}$ is the identity operator. For this case, the tightest estimation of $\mathcal{H}$ is $\mathcal{I}$, which is block diagonal.

Since the coupled function $\phi(u,v)$ consists of two block variables $u$ and $v$, the operators $\mathcal{Q}$ and $\mathcal{W}$ can be decomposed accordingly as $\mathcal{Q} = \begin{pmatrix} \mathcal{Q}_{11} & \mathcal{Q}_{12}\\ \mathcal{Q}_{12}^* & \mathcal{Q}_{22}\end{pmatrix}$ and $\mathcal{W} = \begin{pmatrix} \mathcal{W}_{11} & \mathcal{W}_{12}\\ \mathcal{W}_{12}^* & \mathcal{W}_{22}\end{pmatrix}$, where $\mathcal{W}_{11},~\mathcal{Q}_{11}: \mathcal{U}\to\mathcal{U}$ and $\mathcal{W}_{22}, ~\mathcal{Q}_{22}: \mathcal{V}\to\mathcal{V}$ are self-adjoint positive semidefinite linear operators, and $\mathcal{W}_{12}, ~\mathcal{Q}_{12}:\mathcal{V}\to\mathcal{U}$ are two linear mappings whose adjoints are given by $\mathcal{W}_{12}^*$ and $\mathcal{Q}_{12}^*$, respectively.
Denote $\eta \in [0,1]$ as a constant that satisfies
\begin{equation}\label{cross_control}
|\langle u, (\mathcal{W}_{12}-\mathcal{Q}_{12}) v\rangle |\leq \frac{\eta}{2}(\|u\|^2_{\mathcal{D}_1} + \|v\|^2_{\mathcal{D}_2})\quad\forall\,\mathcal{W}\in \partial^2 \phi(u,v), ~u\in\mathcal{U}, ~v\in\mathcal{V}.
\end{equation}
Note that (\ref{cross_control}) always holds true for $\eta = 1$ according to the Cauchy-Schwarz inequality.

In order to prove the convergence of the proposed majorized ADMM, the following constraint qualification is needed:
\begin{assumption}\label{cq}
There exists $(\hat{u},\hat{v})\in\textup{ri}\;(\textup{dom}(p)\times\textup{dom}(q))$ such that $\mathcal{A}^*\hat{u} + \mathcal{B}^*\hat{v} = c$.
\end{assumption}
Let $\partial p$ and $\partial q$ be the subdifferential mappings of $p$ and $q$, respectively. Define the set-valued mapping $\mathcal{F}$ by
\begin{equation*}
\mathcal{F}(u,v,x): = \nabla \phi(w) + \left(\begin{array}{cc}\partial p(u) + \mathcal{A}x\\[5pt] \partial q(v) + \mathcal{B}x\end{array}\right), \quad (u,v,x)\in\mathcal{U}\times\mathcal{V}\times\mathcal{X}.
\end{equation*}
Under Assumption \ref{cq},  $(\bar{u}, \bar{v})$ is optimal to (\ref{opt}) if and only if there exists $\bar{x}\in\mathcal{X}$ such that the following Karush-Kuhn-Tucker (KKT) condition holds:
\begin{equation}\label{kkt}\left\{\begin{array}{ll}
0\in F(\bar{u},\bar{v},\bar{x}),\\[8pt]
\mathcal{A}^*\bar{u} + \mathcal{B}^*\bar{v} = c,
\end{array}\right.
\end{equation}
which is equivalent to the following variational inequality:
\begin{equation}\label{VI}\begin{array}{rr}
 &(p({u}) + q({v}) ) -( p(\bar u) + q(\bar v) ) + \langle {w}-\bar{w},\nabla\phi(\bar{w})\rangle+ \langle u - \bar{u}, \mathcal{A}\bar{x}\rangle  + \langle v - \bar v , \mathcal{B}\bar{x}\rangle\\[8pt]
 &
 - \langle x - \bar x, \mathcal{A}^*\bar u + \mathcal{B}^*\bar v - c\rangle\geq 0
 \quad\forall (u,v,x)\in \mathcal{U}\times\mathcal{V}\times\mathcal{X}.
 \end{array}
\end{equation}
Motivated by Nesterov's definition of an $\varepsilon$-approximation solution based on the first order optimality condition \cite[Definition 1]{Nesterov_2013} , we say that $(\tilde{u}, \tilde{v}, \tilde{x} ̃)\in\mathcal{U}\times\mathcal{V}\times\mathcal{X}$ is an $\varepsilon$-approximation solution to problem (\ref{opt}) if
\begin{equation}\label{def_approxiamteVI}\begin{array}{ll}
 & (p(\tilde{u}) + q(\tilde{v}) ) -( p(u) + q(v) )+\langle \tilde{w}-w,\nabla\phi(w)\rangle+ \langle   \tilde{u} - u, \mathcal{A}x\rangle  + \langle  \tilde{v} - v , \mathcal{B}x\rangle \\[8pt]
 &
 - \langle \tilde x - x, \mathcal{A}^*u + \mathcal{B}^*v - c\rangle\leq \varepsilon \quad\forall (u,v,x)\in{{B}}(\tilde{u}, \tilde{v}, \tilde{x} ̃),
 \end{array}
\end{equation}
where ${
{B}(\tilde{u}, \tilde{v}, \tilde{x} ̃)}= \{(u,v,x)\in\mathcal{U}\times\mathcal{V}\times\mathcal{X}|\|(u,v,x) - (\tilde u, \tilde v, \tilde x)\| \leq 1\}
$.

Furthermore, since $p$ and $q$ are convex functions, $\partial p(\cdot)$ and $\partial q(\cdot)$ are maximal monotone operators. Then, for any $u,\hat{u}\in \textup{dom}(p)$, $\xi\in\partial p(u)$, and $\hat{\xi}\in\partial p(\hat{u})$, we have
\begin{equation}\label{p_convex}\begin{array}{ll}
&\langle u - \hat{u},\xi - \hat{\xi}\rangle \geq 0,
\end{array}
\end{equation}
and similarly for any $v,\hat{v} \in\textup{dom}(q)$, $\zeta\in\partial q(v)$, and $\hat{\zeta}\in\partial q(\hat{v})$, we have
\begin{equation}\label{q_convex}\begin{array}{ll}
&\langle v - \hat{v}, \zeta- \hat{\zeta}\rangle \geq  0.
\end{array}
\end{equation}
%
\section{A majorized ADMM with coupled objective functions}
In this section, we will first present the framework of our majorized ADMM and then  prove two important inequalities   that play an essential role for our convergence analysis.

Let $\sigma>0$.
For  given $w' = (u',v')\in\mathcal{U}\times \mathcal{V}$,
define the following majorized augmented Lagrangian function associated with (\ref{opt}):
\begin{equation*}\begin{array}{ll}
\widehat{\mathcal{L}}_\sigma(w; (x,w')): =& p(u) + q(v) +\hat{\phi}(w; w') + \langle x,\mathcal{A}^*u + \mathcal{B}^*v - c\rangle +\displaystyle\frac{\sigma}{2}\|\mathcal{A}^*u + \mathcal{B}^*v - c\|^2,
\end{array}
\end{equation*}
where $(w,x) = (u,v,x)\in\mathcal{U}\times\mathcal{V}\times\mathcal{X}$ and
the majorized function $\hat\phi$ is given by (\ref{phi_majorize}).
Then our proposed algorithm works as follows:

\bigskip
\centerline{\fbox{\parbox{\textwidth}{
\textbf{Majorized ADMM:
A majorized ADMM with coupled objective functions} \\[3mm]
Choose an initial point $(u^0, v^0,  x^0)\in\text{dom}(p)\times\text{dom}(q)\times\mathcal{X}$ and parameters $\tau >0$. Let $\mathcal{S}$ and $\mathcal{T}$ be given self-adjoint positive semidefinite linear operators.   Set $k:=0.$ Iterate until convergence:\\[3mm]
\textbf{Step 1.} Compute $u^{k+1} = \displaystyle\arg\min_{u\in\mathcal{U}}\{ \widehat{\mathcal{L}}_\sigma(u,v^k; (x^k,w^k)) + \frac{1}{2}\|u-u^k\|^2_\mathcal{S}\} .$   \\ [8pt]
\textbf{Step 2.} Compute $v^{k+1} = \displaystyle\arg\min_{v\in\mathcal{V}}\{\widehat{\mathcal{L}}_\sigma(u^{k+1},v;(x^k,w^k)) + \frac{1}{2}\|v-v^k\|^2_\mathcal{T}\}$.\\ [8pt]
\textbf{Step 3.} Compute $x^{k+1}=x^k + \tau\sigma(\mathcal{A}^*u^{k+1}+\mathcal{B}^*v^{k+1}-c).$
}}}
\bigskip

In order to simplify  subsequent discussions,  for $k = 0,1,2,\cdots$, define
\begin{equation}\label{notation1}\left\{
\begin{array}{ll}
\tilde{x}^{k+1} := x^k + \sigma(\mathcal{A}^*u^{k+1} + \mathcal{B}^*v^{k+1} - c),\\[8pt]
\Xi_{k+1}  :=   \|v^{k+1} - v^k\|^2_{\mathcal{D}_2+ \mathcal{T}}
+ \eta\|u^{k+1} - u^k\|^2_{\mathcal{D}_1},\\[8pt]
\Theta_{k+1} :=   \|u^{k+1} - u^k\|^2_{\mathcal{S} } + \|v^{k+1} - v^k\|^2_{\mathcal{T}}  + \frac{1}{4}\|w^{k+1}-w^k\|^2_{\mathcal{Q}},\\[8pt]
\Gamma_{k+1} :=  \Theta_{k+1}+  \min(\tau, 1+\tau - \tau^2)\|v^{k+1} - v^k\|^2_{\sigma\mathcal{B}\mathcal{B}^* }  - \|u^{k+1} - u^k\|^2_{\eta\mathcal{D}_1} - \|v^{k+1} - v^k\|^2_{\eta\mathcal{D}_2 }
\end{array}\right.
\end{equation}
and denote for $(u,v,x)\in\mathcal{U}\times\mathcal{V}\times\mathcal{X}$,
\begin{equation}\label{notation2}\left\{\begin{array}{ccl}
\Phi_{k}(u,v,x)& := &(\tau\sigma)^{-1} \|x^{k}-x\|^2 +  \|u^{k}-u\|^2_{\mathcal{D}_1+ \mathcal{S}} + \|v^{k}-v\|^2_{\mathcal{Q}_{22}+ \mathcal{D}_2 + \mathcal{T}}  + \frac{1}{2}\|w^{k}-w\|^2_\mathcal{Q}\\[8pt]
&&
 + \sigma\|\mathcal{A}^*u+\mathcal{B}^*v^k-c\|^2,\\[8pt]
\Psi_{k}(u,v,x) &:= & \Phi_{k}(u,v,x)+\|w^{k}-w\|^2_{\mathcal{Q}}+  \max(1-\tau,1-\tau^{-1})\sigma\|\mathcal{A}^*u^{k}   + \mathcal{B}^*v^{k} - c\|^2.
\end{array}\right.
\end{equation}

\begin{prop}\label{decrease}
 Suppose that the solution set of problem (\ref{opt}) is nonempty and Assumption \ref{cq} holds. Assume that $\mathcal{S}$ and $\mathcal{T}$ are chosen such that the sequence $\{(u^k,v^k,x^k)\}$ is well defined. Then the following conclusions hold:\\[8pt]
 (i) For $\tau\in(0,1]$,  we have that for any $k\geq 0$ and $(u,v,x)\in\mathcal{U}\times\mathcal{V}\times\mathcal{X}$,
 \begin{equation}\label{decrease_11}
\begin{array}{ll}
&\disp(p(u^{k+1}) + q(v^{k+1}) ) -(p(u) + q(v) ) +\langle w^{k+1}-w,\nabla\phi(w)\rangle+ \langle   u^{k+1} - u, \mathcal{A}x\rangle  + \langle  v^{k+1} - v , \mathcal{B}x\rangle\\[8pt]
&\disp
 - \langle  \tilde{x}^{k+1} - x , \mathcal{A}^*u + \mathcal{B}^*v - c\rangle +
\frac{1}{2}(\Phi_{k+1}({u}, {v}, {x}) - \Phi_k({u}, {v}, {x}))\\[8pt]

\leq  & \disp-\frac{1}{2}(\Theta_{k+1}+  \sigma\|\mathcal{A}^*u^{k+1} + \mathcal{B}^*v^k - c\|^2
 +(1-\tau)\sigma\|\mathcal{A}^*u^{k+1} + \mathcal{B}^*v^{k+1} - c\|^2  ).
\end{array}
\end{equation}
(ii) For $\tau \geq 0$, we have that for any $k\geq 1$ and $(u,v,x)\in\mathcal{U}\times\mathcal{V}\times\mathcal{X}$,
\begin{equation}\label{decrease_22}
\begin{array}{ll}
&(p(u^{k+1}) + q(v^{k+1}) ) -(p(u) + q(v) ) +\langle w^{k+1}-w,\nabla\phi(w)\rangle+ \langle   u^{k+1} - u, \mathcal{A}x\rangle  + \langle  v^{k+1} - v , \mathcal{B}x\rangle\\[8pt]
&\disp
 - \langle  \tilde{x}^{k+1} - x , \mathcal{A}^*u + \mathcal{B}^*v - c\rangle +
\frac{1}{2}(\Psi_{k+1}({u}, {v}, {x}) + \Xi_{k+1} - (\Psi_k({u}, {v}, {x}) +\Xi_k))\\[8pt]
\leq &\disp
 -\frac{1}{2}(\Gamma_{k+1}  +\min(1,1+\tau^{-1} - \tau)\sigma\|\mathcal{A}^*u^{k+1}   + \mathcal{B}^*v^{k+1} - c\|^2).
\end{array}
\end{equation}
\end{prop}
{\bf Proof.}
 In the majorized  ADMM  iteration scheme, the optimality condition for $(u^{k+1},v^{k+1})$ is
\begin{equation}\label{opt_1}\left\{\begin{array}{ll}
0\in&\partial p(u^{k+1})  + \nabla_u \phi(w^k)  + \mathcal{A}x^k + \sigma\mathcal{A}(\mathcal{A}^*u^{k+1} + \mathcal{B}^*v^k - c) + (\mathcal{Q}_{11} + \mathcal{D}_1+ \mathcal{S})(u^{k+1} - u^k),\\[8pt]
0\in&\partial q(v^{k+1})+   \nabla_v \phi(w^k) + \mathcal{B}x^k + \sigma\mathcal{B}(\mathcal{A}^*u^{k+1} + \mathcal{B}^*v^{k+1} - c) + (\mathcal{Q}_{22} + \mathcal{D}_2+ \mathcal{T})(v^{k+1} - v^k) \\[8pt]
&
+ \mathcal{Q}_{12}^*(u^{k+1} - u^k).
\end{array}\right.
\end{equation}
Denote
\begin{equation*}\left\{\begin{array}{ll}
a^{k+1} &= -x^{k+1} - (1-\tau)\sigma(\mathcal{A}^*u^{k+1}+\mathcal{B}^*v^{k+1}-c) - \sigma\mathcal{B}^*(v^k - v^{k+1}),\\[8pt]
b^{k+1} & = -x^{k+1} - (1-\tau)\sigma (\mathcal{A}^*u^{k+1}+\mathcal{B}^*v^{k+1}-c).
\end{array}\right.
\end{equation*}
Then by noting that 
$$x^{k+1}= x^k + \tau\sigma(\mathcal{A}^*u^{k+1}+\mathcal{B}^*v^{k+1}-c),$$
we can rewrite (\ref{opt_1}) as 
\begin{equation}\label{opt_2}\left\{\begin{array}{ll}
&\mathcal{A}a^{k+1}  -\nabla_u \phi(w^k) - (\mathcal{Q}_{11} + \mathcal{D}_1+ \mathcal{S})(u^{k+1} - u^k)\in\partial p(u^{k+1}),\\[8pt]

&\mathcal{B}b^{k+1} - \nabla_v \phi(w^k)-(\mathcal{Q}_{22} + \mathcal{D}_2 +  \mathcal{T})(v^{k+1} - v^k) - \mathcal{Q}_{12}^*(u^{k+1} - u^k)  \in\partial q(v^{k+1}).
\end{array}\right.
\end{equation}
Therefore, by the convexity of $p$ and $q$, we have that for any $u\in\mathcal{U}$ and $v\in\mathcal{V}$,
\begin{equation}\label{pqineq}\left\{\begin{array}{l}
p(u) \geq p(u^{k+1}) + \langle u - u^{k+1},\mathcal{A}a^{k+1}  -\nabla_u \phi(w^k) - (\mathcal{Q}_{11}+\mathcal{D}_1+ \mathcal{S})(u^{k+1} - u^k)\rangle,\\[8pt]

q(v) \geq q(v^{k+1}) + \langle v - v^{k+1}, \mathcal{B}b^{k+1}  - \nabla_v \phi(w^k)-(\mathcal{Q}_{22} + \mathcal{D}_2 +  \mathcal{T})( v^{k+1}-v^k)-\mathcal{Q}^*_{12}(u^{k+1}- u^{k})\rangle.
\end{array}\right.
\end{equation}
By noting the relationship between $x^{k+1}, \tilde{x}^{k+1}$ and $x^k$, we obtain from  the above inequalities that for any $(u,v)\in\mathcal{U}\times\mathcal{V}$,
\begin{equation}\label{main_ineq}\begin{array}{ll}
&(p(u^{k+1}) + q(v^{k+1}) ) -(p(u) + q(v) ) +\langle w^{k+1}-w,\nabla\phi(w)\rangle+ \langle   u^{k+1} - u, \mathcal{A}x\rangle  + \langle  v^{k+1} - v , \mathcal{B}x\rangle\\[8pt]
&
 - \langle  \tilde{x}^{k+1} - x , \mathcal{A}^*u + \mathcal{B}^*v - c\rangle \\[8pt]
\leq & \sigma\langle \mathcal{B}^*(v^{k+1} - v^k),\mathcal{A}^*(u^{k+1}-u)\rangle

-\langle w^{k+1} - w, \nabla\phi(w^{k}) - \nabla\phi(w)\rangle  - \langle \mathcal{Q}_{12}^*(u^{k+1} - u^k), v^{k+1}-v\rangle\\[8pt]
&
-\langle (\mathcal{Q}_{11}+\mathcal{D}_1 + \mathcal{S})(u^{k+1}-u^k),u^{k+1} - u\rangle
- \langle (\mathcal{Q}_{22}+\mathcal{D}_2 + \mathcal{T})(v^{k+1} - v^k), v^{k+1}-v\rangle\\[8pt]
&
 - (\tau\sigma)^{-1}\langle x^{k+1} - x^k, x^{k+1} - x\rangle - (1-\tau)\sigma\|\mathcal{A}^*u^{k+1} + \mathcal{B}^*v^{k+1} - c\|^2.
\end{array}
\end{equation}
By taking $(w, w') = (\bar{w}, w^k)$ and $(w^{k+1}, \bar{w})$ in (\ref{phi_convex}),  we  know that
\begin{equation*}\begin{array}{rll}
\phi({w}) &\geq &  \phi(w^k) + \langle \nabla \phi(w^k), {w} - w^k\rangle + \frac{1}{2}\|{w} - w^k\|^2_{\mathcal{Q}}  ,\\ [8pt]
\phi(w^{k+1}) & \geq & \phi({w}) + \langle \nabla \phi({w}),w^{k+1} - {w}\rangle + \frac{1}{2}\|w^{k+1} - {w}\|^2_{\mathcal{Q}} .
\end{array}
\end{equation*}
By taking $(w, w') = (w^{k+1},  w^k)$ in (\ref{phi_majorize}), we can also get that
\begin{equation*}
\phi(w^{k+1}) \leq \phi(w^{k})  + \langle \nabla \phi(w^{k}), w^{k+1} - w^k\rangle +\frac{1}{2}\|w^{k+1} - w^k\|^2_{\mathcal{Q} +\mathcal{H}}.
\end{equation*}
Putting the above three inequalities together, we get
\begin{equation}\label{fg_estimate}
 \langle \nabla \phi(w^k) - \nabla \phi({w}), w^{k+1}-w\rangle
  \geq
\frac{1}{2}(\|w^k-w\|^2_{\mathcal{Q}}+\|w^{k+1}-w\|^2_{\mathcal{Q}}) -\frac{1}{2}\|w^{k+1} - w^k\|^2_{\mathcal{Q} +\mathcal{H}}.
\end{equation}
Note that
\begin{equation}\label{equiv}
\begin{array}{ll}
&\disp\frac{1}{2}(\|w^{k+1} - w^k\|^2_{\mathcal{Q}} - \|w^{k+1}-w\|^2_\mathcal{Q}  - \|w^k-w\|^2_\mathcal{Q}) - \langle \mathcal{Q}_{11}(u^{k+1} - u^k), u^{k+1}-u\rangle\\[8pt]
& - \langle \mathcal{Q}_{22}(v^{k+1} - v^k), v^{k+1}-v\rangle
 - \langle \mathcal{Q}^*_{12}( u^{k+1} - u^k), v^{k+1}-v\rangle\\[8pt]
 = &\disp \frac{1}{2}(\|w^{k+1} - w^k\|^2_{\mathcal{Q}} - \|w^{k+1}-w\|^2_{\mathcal{Q}}  - \|w^k-w\|^2_{\mathcal{Q}})  - \langle w^{k+1} - w^k,\mathcal{Q}(w^{k+1}-w)\rangle \\[8pt]
 &
 + \langle \mathcal{Q}_{12}(v^{k+1} - v^k), u^{k+1}-u\rangle\\[8pt]
  = &-\|w^{k+1}-w\|^2_{\mathcal{Q}}+ \langle \mathcal{Q}_{12}(v^{k+1} - v^k), u^{k+1}-u\rangle.
  \end{array}
\end{equation}
Substituting (\ref{fg_estimate}) and  (\ref{equiv}) into (\ref{main_ineq}) and by the assumption (\ref{diag_H}), we can further obtain that
\begin{equation}\label{main_ineq_2}\begin{array}{ll}
&(p(u^{k+1}) + q(v^{k+1}) ) -(p(u) + q(v) ) +\langle w^{k+1}-w,\nabla\phi(w)\rangle+ \langle   u^{k+1} - u, \mathcal{A}x\rangle  + \langle  v^{k+1} - v , \mathcal{B}x\rangle\\[8pt]
&
 - \langle  \tilde{x}^{k+1} - x , \mathcal{A}^*u + \mathcal{B}^*v - c\rangle \\[8pt]
\leq
&
\sigma\langle \mathcal{B}^*(v^{k+1} - v^k),\mathcal{A}^*(u^{k+1}-u)\rangle- \|w^{k+1}-w\|^2_{\mathcal{Q}}+\frac{1}{2}\|w^{k+1} - w^k\|^2_{\text{Diag}\,(\mathcal{D}_1, \mathcal{D}_2)}\\[8pt]
&
+\langle v^{k+1}-v^k, \mathcal{Q}_{12}^*(u^{k+1} - u)\rangle
-\langle (\mathcal{D}_1 + \mathcal{S})(u^{k+1}-u^k),u^{k+1} - u\rangle\\[8pt]
&
- \langle (\mathcal{D}_2 + \mathcal{T})(v^{k+1} - v^k), v^{k+1}-v\rangle
 - (\tau\sigma)^{-1}\langle x^{k+1} - x^k, x^{k+1} - x\rangle\\[8pt]
 & - (1-\tau)\sigma\|\mathcal{A}^*u^{k+1} + \mathcal{B}^*v^{k+1} - c\|^2.
\end{array}
\end{equation}
Recall that for any $\xi$, $\zeta$ in the same space and a self-adjoint positive semidefinite operator $\mathcal{G}$, it always holds that
\begin{equation}\label{inner_product_eq}
\langle \xi, \mathcal{G}\zeta\rangle  = \frac{1}{2}(\|\xi\|_\mathcal{G}^2 + \|\zeta\|_\mathcal{G}^2 - \|\xi-\zeta\|_\mathcal{G}^2) = \frac{1}{2}(\|\xi+\zeta\|_\mathcal{G}^2 - \|\xi\|_\mathcal{G}^2 - \|\zeta\|_\mathcal{G}^2).
\end{equation}
Then we can get that
\begin{equation}\label{cross_st}\begin{array}{rcl}
\langle (\mathcal{D}_1 + \mathcal{S})(u^{k+1} - u^k), u^{k+1} - u\rangle &=& \displaystyle\frac{1}{2}(\|u^{k+1} - u^k\|_{\mathcal{D}_1+\mathcal{S}}^2 + \|u^{k+1} - u\|_{\mathcal{D}_1+\mathcal{S}}^2 - \|u^k - u\|_{\mathcal{D}_1+\mathcal{S}}^2),\\[8pt]
\langle (\mathcal{D}_2 + \mathcal{T})(v^{k+1} - v^k), v^{k+1} - v\rangle &=& \displaystyle\frac{1}{2}(\|v^{k+1} - v^k\|_{\mathcal{D}_2 + \mathcal{T}}^2 + \|v^{k+1} - v\|_{\mathcal{D}_2 + \mathcal{T}}^2 - \|v^k - v\|^2_{\mathcal{D}_2 + \mathcal{T}}),\\[8pt]
\langle x^{k+1} - x^k, x^{k+1} - x\rangle &=& \displaystyle \frac{1}{2}(\|x^{k+1} - x^k\|^2 + \|x^{k+1} - x\|^2 - \|x^k - x\|^2),\\[8pt]
\langle \mathcal{Q}_{22}(v^{k+1} - v^{k}), v^{k+1}-v  \rangle & = &\displaystyle \frac{1}{2}(\|v^{k+1} - v^k\|_{\mathcal{Q}_{22}}^2 + \|v^{k+1} - v\|_{\mathcal{Q}_{22}}^2 - \|v^k - v\|_{\mathcal{Q}_{22}}^2).
\end{array}
\end{equation}
{\bf (i)} Assume that $\tau\in(0,1]$. 
By using the last equation in (\ref{cross_st}), we can obtain that
\begin{equation}\label{cross_1}\begin{array}{lll}
\inprod{v^{k+1} - v^{k}}{\mathcal{Q}_{12}^*(u^{k+1}-u) }
&= & \Inprod{ \left(\begin{array}{c} 0\\ v^{k+1} - v^{k}\end{array}\right)}{\mathcal{Q}( w^{k+1}-w )}- \langle \mathcal{Q}_{22}(v^{k+1} - v^{k}), v^{k+1}-v  \rangle \\[8pt]

&\leq & \disp\frac{1}{2}(\|v^{k+1} - v^k\|^2_{\mathcal{Q}_{22}} + \|w^{k+1}-w\|^2_{\mathcal{Q}})  - \frac{1}{2}(\|v^{k+1} - v^k\|^2_{\mathcal{Q}_{22}} \\[8pt]
&&
+ \|v^{k+1}-v\|^2_{\mathcal{Q}_{22}} - \|v^k-v\|^2_{\mathcal{Q}_{22}})\\[8pt]

&= & \disp\frac{1}{2}\|w^{k+1}-w\|^2_{\mathcal{Q}}  + \frac{1}{2}(\|v^{k}-v\|^2_{\mathcal{Q}_{22}} - \|v^{k+1}-v\|^2_{\mathcal{Q}_{22}}),
\end{array}
\end{equation}
where the inequality is obtained by the Cauchy-Schwarz inequality.
By some simple manipulations  we can also see that
\begin{equation}\label{cross_2}\begin{array}{ll}
\sigma\langle \mathcal{B}^*(v^{k+1} - v^k), \mathcal{A}^*(u^{k+1}-u)\rangle
& = \disp\frac{\sigma}{2}( \|\mathcal{A}^*u^{k+1} + \mathcal{B}^*v^{k+1} - c\|^2 - \|\mathcal{A}^*u^{k+1} + \mathcal{B}^*v^k - c\|^2)\\[8pt]
& \disp\quad +\frac{\sigma}{2}(\|\mathcal{A}^*u + \mathcal{B}^*v^k-c\|^2 -\|\mathcal{A}^*u + \mathcal{B}^*v^{k+1}-c\|^2).
\end{array}
\end{equation}
Finally, by substituting (\ref{cross_st}),   (\ref{cross_1}) and (\ref{cross_2}) into (\ref{main_ineq_2}) and recalling the definition of $\Phi_{k+1}(\cdot, \cdot, \cdot)$ and $\Theta_{k+1}$ in (\ref{notation1}) and (\ref{notation2}), we have that
\begin{equation*}\begin{array}{ll}
&(p(u^{k+1}) + q(v^{k+1}) ) -(p(u) + q(v) ) +\langle w^{k+1}-w,\nabla\phi(w)\rangle+ \langle   u^{k+1} - u, \mathcal{A}x\rangle  + \langle  v^{k+1} - v , \mathcal{B}x\rangle\\[8pt]
&
 \disp - \langle  \tilde{x}^{k+1} - x , \mathcal{A}^*u + \mathcal{B}^*v - c\rangle +
\frac{1}{2}(\Phi_{k+1}({u}, {v}, {x}) - \Phi_k({u}, {v}, {x}))\\[8pt]
&  \leq  \disp - \frac{1}{2} ( \|u^{k+1} - u^k\|^2_{\mathcal{S} } + \|v^{k+1} - v^k\|^2_{\mathcal{T}}  + \frac{1}{2}\|w^{k+1}-w\|^2_{\mathcal{Q}}+\frac{1}{2} \|w^k-w\|^2_{\mathcal{Q}}
+  \sigma\|\mathcal{A}^*u^{k+1} + \mathcal{B}^*v^k - c\|^2 \\[8pt]
&
\quad\quad\quad+(1-\tau)\sigma\|\mathcal{A}^*u^{k+1} + \mathcal{B}^*v^{k+1} - c\|^2)\\[8pt]
& \leq \disp  -  \frac{1}{2}(\Theta_{k+1}  +  \sigma\|\mathcal{A}^*u^{k+1} + \mathcal{B}^*v^k - c\|^2+ (1-\tau)\sigma\|\mathcal{A}^*u^{k+1} + \mathcal{B}^*v^{k+1} - c\|^2),
\end{array}
\end{equation*}
where the last inequality comes from the fact that
\begin{equation*}
\frac{1}{2}\|w^{k+1}-w\|^2_{\mathcal{Q}} + \frac{1}{2}\|w^{k}-w\|_{\mathcal{Q}}^2\geq \frac{1}{4}\|w^{k+1} - w^k\|^2_{\mathcal{Q}}.
\end{equation*}
This completes the proof of part (i).\\[8pt]
{\bf (ii)} Assume that $\tau\geq 0$. In this part, we first reformulate (\ref{main_ineq_2}) as
\begin{equation}\label{main_ineq_3}\begin{array}{ll}
&(p(u^{k+1}) + q(v^{k+1}) ) -(p(u) + q(v) ) +\langle w^{k+1}-w,\nabla\phi(w)\rangle+ \langle   u^{k+1} - u, \mathcal{A}x\rangle  + \langle  v^{k+1} - v , \mathcal{B}x\rangle\\[8pt]
&
 - \langle  \tilde{x}^{k+1} - x , \mathcal{A}^*u + \mathcal{B}^*v - c\rangle \\[8pt]
\leq
&
\disp\sigma\langle \mathcal{B}^*(v^{k+1} - v^k),\mathcal{A}^*u^{k+1} + \mathcal{B}^*v^{k+1} - c\rangle+\frac{\sigma}{2}(\|\mathcal{A}^*u + \mathcal{B}^*v^k-c\|^2 -\|\mathcal{A}^*u + \mathcal{B}^*v^{k+1}-c\|^2) \\[8pt]
&
\disp- \frac{1}{2}\|v^{k+1} - v^k\|^2_{\sigma\mathcal{B}\mathcal{B}^*}

- \|w^{k+1}-w\|^2_{\mathcal{Q}}+\frac{1}{2}\|w^{k+1} - w^k\|^2_{\text{Diag}\,(\mathcal{D}_1, \mathcal{D}_2)}  -   \langle \mathcal{Q}_{22}(v^{k+1} - v^k), v^{k+1}-v\rangle\\[8pt]
&
+ \langle v^{k+1}-v^k, \mathcal{Q}_{12}^*(u^{k+1} - u) + \mathcal{Q}_{22}(v^{k+1} - v)\rangle
-\langle (\mathcal{D}_1 + \mathcal{S})(u^{k+1}-u^k),u^{k+1} - u\rangle\\[8pt]
&

- \langle (\mathcal{D}_2 + \mathcal{T})(v^{k+1} - v^k), v^{k+1}-v\rangle
 - (\tau\sigma)^{-1}\langle x^{k+1} - x^k,  x^{k+1} - x\rangle \\[8pt]
&
- (1-\tau)\sigma\|\mathcal{A}^*u^{k+1} + \mathcal{B}^*v^{k+1} - c\|^2.
\end{array}
\end{equation}
Next we shall estimate the following cross term 
$$\sigma\langle \mathcal{B}^*(v^{k+1} - v^k), \mathcal{A}^*u^{k+1}  + \mathcal{B}^*v^{k+1} - c\rangle  + \langle  v^{k+1} - v^k, \mathcal{Q}_{12}^*(u^{k+1}-u) + \mathcal{Q}_{22}(v^{k+1}-v)\rangle.$$
It follows from  (\ref{opt_2}) that
\begin{equation}\label{second_opt}
\left\{\begin{array}{ll}
\mathcal{B}b^{k+1} - \nabla_v \phi(w^k)-(\mathcal{Q}_{22}+\mathcal{D}_2 +  \mathcal{T})(v^{k+1} - v^k) - \mathcal{Q}_{12}^*(u^{k+1} - u^k)  \in\partial q(v^{k+1}), \\[8pt]
\mathcal{B}b^{k} - \nabla_v \phi(w^{k-1})-(\mathcal{Q}_{22}+\mathcal{D}_2 +  \mathcal{T})(v^{k} - v^{k-1})-\mathcal{Q}_{12}^*(u^k-u^{k-1})  \in\partial q(v^{k}).

\end{array}\right.
\end{equation}
Since $\nabla \phi$ is globally Lipschitz continuous, it is known from Clarke's Mean-Value Theorem \cite[Proposition 2.6.5]{Clarke}  that there exists a self-adjoint and positive semidefinite operator \\
$\mathcal{W}^k \in \text{conv}\{\partial^2\phi([w^{k-1}, w^k])\} $ such that
\begin{equation*}
\nabla\phi(w^k) - \nabla\phi(w^{k-1})= \mathcal{W}^k(w^k-w^{k-1}),
\end{equation*}
{where the set $\text{conv}\{\partial^2\phi[w^{k-1},w^k]\}$ denotes the convex hull of all points 
$\cW\in \partial^2\phi(z)$ for any $z\in [w^{k-1},w^k]$.}
Denote $\mathcal{W}^k: = \begin{pmatrix} \mathcal{W}_{11}^k & \mathcal{W}_{12}^k \\ (\mathcal{W}_{12}^k)^* & \mathcal{W}_{22}^k\end{pmatrix}$, where $\mathcal{W}_{11}^k: \mathcal{U}\to\mathcal{U}$, $\mathcal{W}_{22}^k: \mathcal{V}\to\mathcal{V}$ are self-adjoint positive semidefinite operators and $\mathcal{W}_{12}^k: \mathcal{U}\to\mathcal{V}$ is a linear operator. {Substituting (\ref{second_opt}) into (\ref{q_convex}) at $ v = v^{k+1}$ and $\hat{v} = v^k$, we obtain that
\begin{equation*}\begin{array}{ll}
&\langle \mathcal{B}(b^{k+1} - b^k), v^{k+1} - v^k\rangle -\langle \mathcal{Q}_{22}(v^{k+1} - v^{k}) +\mathcal{Q}_{12}^*  (u^{k+1} - u^{k}), v^{k+1} - v^k\rangle \\[8pt]

\geq &\langle \nabla_v\phi(w^k)  - \nabla_v \phi(w^{k-1}), v^{k+1} -v^k\rangle -\langle (\mathcal{Q}_{22} + \mathcal{D}_2+  \mathcal{T})(v^{k} - v^{k-1}), v^{k+1} - v^k\rangle \\[8pt]
&
\quad + \|v^{k+1} - v^k\|^2_{\mathcal{T}+\mathcal{D}_2}  -\langle u^k-u^{k-1}, \mathcal{Q}_{12}(v^{k+1} - v^k)\rangle\\[8pt]
= & \langle u^k-u^{k-1},(\mathcal{W}_{12}^k-\mathcal{Q}_{12})( v^{k+1} - v^k)\rangle -\langle ( \mathcal{Q}_{22}+ \mathcal{D}_2+\mathcal{T}  - \mathcal{W}_{22}^k)(v^{k} - v^{k-1}), v^{k+1} - v^k\rangle\\[8pt]
&
\quad + \|v^{k+1} - v^k\|^2_{\mathcal{T}+\mathcal{D}_2} \\[8pt]

 \geq & \disp-\frac{\eta}{2}(\|u^k - u^{k-1}\|^2_{\mathcal{D}_1} + \|v^{k+1} - v^k\|^2_{\mathcal{D}_2}) - \frac{1}{2}(\|v^{k+1} - v^k\|^2_{\mathcal{T}+ \mathcal{D}_2}
  + \|v^k-v^{k-1}\|^2_{\mathcal{T}+ \mathcal{D}_2})\\[8pt]
  &
 \quad+\|v^{k+1} - v^k\|^2_{\mathcal{T}+\mathcal{D}_2}\\[8pt]
= & \disp\frac{1}{2}\|v^{k+1} - v^k\|^2_{\mathcal{T} +(1-\eta)\mathcal{D}_2} - \frac{1}{2}\|v^k - v^{k-1}\|^2_{\mathcal{T} +\mathcal{D}_2} - \frac{\eta}{2}\|u^k-u^{k-1}\|^2_{\mathcal{D}_1},
 \end{array}
\end{equation*}
where the second inequality is obtained from (\ref{cross_control}) and the fact that $\mathcal{W}_{22}^k\succeq \mathcal{Q}_{22}$.
Therefore, with $\mu_{k+1} = (1-\tau)\sigma\langle   \mathcal{B}^*(v^{k+1} - v^k), \mathcal{A}^*u^{k}  + \mathcal{B}^*v^{k} - c\rangle$,  the cross term  can be estimated as
\begin{equation}\label{lim_10}\begin{array}{ll}
&\sigma\langle \mathcal{B}^*(v^{k+1} - v^k), \mathcal{A}^*u^{k+1}  + \mathcal{B}^*v^{k+1} - c\rangle + \langle \mathcal{Q}_{12}^*(u^{k+1}-u) + \mathcal{Q}_{22}(v^{k+1}-v), v^{k+1} - v^k\rangle\\[8pt]

 =& (1-\tau)\sigma\langle   \mathcal{B}^*(v^{k+1} - v^k), \mathcal{A}^*u^{k}  + \mathcal{B}^*v^{k} - c\rangle - \langle
\mathcal{B}^*(v^{k+1} - v^k),b^{k+1}-b^{k}\rangle\\ [8pt]
& + \langle \mathcal{Q}_{12}^*(u^k-u) + \mathcal{Q}_{22}(v^k-v), v^{k+1} - v^k\rangle
  + \langle \mathcal{Q}_{12}^*(u^{k+1} - u^k) + \mathcal{Q}_{22}(v^{k+1} - v^k), v^{k+1} - v^k\rangle\\[8pt]

\leq &\disp\mu_{k+1} + \frac{1}{2}(\|w^k-w\|^2_{\mathcal{Q}} + \|v^{k+1} - v^k\|^2_{\mathcal{Q}_{22}})
- \frac{1}{2}\|v^{k+1} - v^k\|^2_{\mathcal{T} +(1-\eta)\mathcal{D}_2} + \frac{1}{2}\|v^k - v^{k-1}\|^2_{\mathcal{T} + \mathcal{D}_2}\\[8pt]
&
+\disp \frac{\eta}{2}\|u^k-u^{k-1}\|^2_{\mathcal{D}_1}.

\end{array}
\end{equation}
Finally, by the Cauchy-Schwarz inequality we know that
\begin{equation}\label{lim_11}
\mu_{k+1} \leq\left\{\begin{array}{ll}
\disp \frac{1}{2}(1-\tau)\sigma(\|\mathcal{B}^*(v^{k+1} - v^k)\|^2 + \|\mathcal{A}^*u^k + \mathcal{B}^*v^k - c\|^2), \quad\tau \in(0,1],
\\[8pt]
\disp\frac{1}{2}(\tau-1)\sigma(\tau\|\mathcal{B}^*(v^{k+1} - v^k)\|^2 +\tau^{-1} \|\mathcal{A}^*u^k + \mathcal{B}^*v^k - c\|^2), \quad \tau >1.
\end{array}\right.
\end{equation}
Substituting (\ref{cross_st}), (\ref{lim_10}) and (\ref{lim_11})  into (\ref{main_ineq_3}), we can obtain (\ref{decrease_22}). This completes the proof of part (ii).
\qed

\section{Convergence analysis}
With all the preparations given  in the previous sections, we can now discuss the main convergence results of our paper.

\subsection{The global convergence}
First we prove that under mild conditions, the iteration sequence $\{(u^k, v^k, x^k)\}$ generated by the majorized ADMM with $\tau\in(0,\frac{1 +\sqrt{5}}{2})$ converges to an optimal solution of problem (\ref{opt}) and its dual.

Let $\bar{w} = (\bar{u}, \bar{v})\in\mathcal{U}\times\mathcal{V}$ be an optimal solution of (\ref{opt}) and $\bar{x}\in\mathcal{X}$ be the corresponding optimal multiplier.
For $k = 0,1,2,\cdots$,  define
\begin{equation*}\begin{array}{cc}
u_e^k = u^k - \bar{u},  \quad v_e^k = v^k - \bar{v}, \quad w_e^k = w^k - \bar{w}, \quad  x_e^k = x^k - \bar{x}.
\end{array}
\end{equation*}
\begin{theorem}\label{global_convergence}
Suppose that the solution set of (\ref{opt}) is nonempty and Assumption \ref{cq} holds. Assume that $\mathcal{S}$ and $\mathcal{T}$ are chosen such that \[\mathcal{Q}_{11} + \sigma\mathcal{A}\mathcal{A}^*+\mathcal{S}\succ 0, \quad \mathcal{Q}_{22} + \sigma\mathcal{B}\mathcal{B}^*+\mathcal{T}\succ 0.\]
(i) Assume that $\tau\in(0,1]$. If for any $w = \left(\begin{array}{c}u\\v\end{array}\right) \in\mathcal{U}\times\mathcal{V}$, it holds that
\begin{equation}\label{case1_assumption}
\langle w, [\mathcal{Q}+\text{Diag}\,(\mathcal{S}+(1-\tau)\sigma\mathcal{A}\mathcal{A}^*, \mathcal{T} + (1-\tau)\sigma\mathcal{B}\mathcal{B}^*)]w\rangle = 0 \Rightarrow \|u\|\|v\| = 0,
\end{equation}
then the generated sequence $\{(u^k, v^k)\}$ converges to an optimal solution of (\ref{opt}) and $\{x^k\}$ converges to the corresponding optimal multiplier.\\[8pt]
(ii) Assume that $\tau\in (0, \frac{1+\sqrt{5}}{2})$. Under the conditions that
\begin{equation}\label{convergence_2}\begin{array}{cc}
\disp\mathcal{M}: = \frac{1}{4}\mathcal{Q} + \text{Diag}\,(\mathcal{S}-\eta\mathcal{D}_1, ~\mathcal{T} -\eta\mathcal{D}_2 )\succeq 0,\\[8pt]
\disp\frac{1}{4} \mathcal{Q}_{11}+\mathcal{S} + \sigma\mathcal{A}\mathcal{A}^* - \eta\mathcal{D}_1\succ 0,\quad
 \frac{1}{4}\mathcal{Q}_{22}+\mathcal{T} + \sigma\mathcal{B}\mathcal{B}^* - \eta\mathcal{D}_2\succ 0
\end{array}
\end{equation}
and for any $w = \left(\begin{array}{c}u\\v\end{array}\right) \in\mathcal{U}\times\mathcal{V}$, it holds that
\begin{equation}\label{case2_assumption}
\langle w, [\mathcal{M}+\sigma\text{Diag}\,(\mathcal{A}\mathcal{A}^*, \mathcal{B}\mathcal{B}^* )]w\rangle = 0 \Rightarrow \|u\|\|v\| = 0,
\end{equation}
the generated sequence $\{(u^k, v^k)\}$ converges to an optimal solution of (\ref{opt}) and $\{x^k\}$ converges to the corresponding optimal  multiplier.
\end{theorem}
{\bf Proof.} 
{\bf (i)} Let $\tau\in(0,1]$. By letting $(u,v,x) = (\bar{u}, \bar{v}, \bar{x})$ in inequality (\ref{decrease_11}) and the optimality condition (\ref{VI}),
 we can obtain that for any $k\geq 0$,
 \begin{equation}\label{convergence_1}
 \begin{array}{ll}
&\Phi_{k+1}(\bar{u}, \bar{v}, \bar{x}) - \Phi_k(\bar{u}, \bar{v}, \bar{x})\\[8pt]
\leq  & -(\Theta_{k+1}+  \sigma\|\mathcal{A}^*u^{k+1} + \mathcal{B}^*v^k - c\|^2
 +(1-\tau)\sigma\|\mathcal{A}^*u^{k+1} + \mathcal{B}^*v^{k+1} - c\|^2  ).
\end{array}
\end{equation}
The above inequality shows that $\{\Phi_{k+1}(\bar{u}, \bar{v}, \bar{x})\}$ is bounded,
which implies that $\{\|x^{k+1}\|\}$, $\{\|w_e^{k+1}\|_{\mathcal{Q}}\}$, $\{\|u_e^{k+1}\|_{S}\}$ and $\{\|v_e^{k+1}\|_{\mathcal{Q}_{22} + \sigma\mathcal{B}\mathcal{B}^*+\mathcal{T}  }\}$ are all bounded. From the positive definiteness of $\mathcal{Q}_{22} + \sigma\mathcal{B}\mathcal{B}^*+\mathcal{T}$, we can see that $\{\|v_e^{k+1}\|\}$ is bounded.   By using the inequalities
\begin{equation*}\begin{array}{ll}
\|\mathcal{A}^*u_e^{k+1}\|&\leq \|\mathcal{A}^*u_e^{k+1} + \mathcal{B}^*v_e^{k+1}\| + \|\mathcal{B}^*v_e^{k+1}\|\\[8pt]
& \leq \tau\sigma(\|x_e^{k+1}\| + \|x_e^k\|) + \|\mathcal{B}^*v_e^{k+1}\|,\\[8pt]
\|u^{k+1}_e\|_{\mathcal{Q}_{11}} & \leq \|w_e^{k+1}\|_{\mathcal{Q}} + \|v_e^{k+1}\|_{\mathcal{Q}_{22}},
\end{array}
\end{equation*}
we know that the sequence $\{\|u_e^{k+1}\|_{\sigma\mathcal{A}\mathcal{A}^* + \mathcal{Q}_{11}}\}$ is also bounded. Therefore, $\{\|u_e^{k+1}\|_{\mathcal{Q}_{11} + \sigma\mathcal{A}\mathcal{A}^*+\mathcal{S}} \}$ is bounded. By the positive definiteness of $\mathcal{Q}_{11} + \sigma\mathcal{A}\mathcal{A}^*+\mathcal{S}$, we know that $\{\|u_e^{k+1}\|\}$ is bounded. On the whole, the sequence $\{(u^k,v^k,x^k)\}$ is bounded. Thus, there exists a subsequence $\{(u^{k_i},v^{k_i}, x^{k_i})\}$ converging to a cluster point, say $(u^{\infty},v^{\infty},x^{\infty})$.
Next we will prove that $(u^{\infty},v^{\infty})$ is optimal to (\ref{opt}) and $x^\infty$ is the corresponding optimal multiplier. The inequality (\ref{convergence_1}) also implies that $$\displaystyle \lim_{k\to\infty} (\Theta_{k+1} +  \sigma\|\mathcal{A}^*u^{k+1} + \mathcal{B}^*v^k - c\|^2
+(1-\tau)\sigma\|\mathcal{A}^*u^{k+1} + \mathcal{B}^*v^{k+1} - c\|^2)  = 0,$$ which is equivalent to
\begin{equation}\label{lim_1}\begin{array}{ll}
\displaystyle \lim_{k\to\infty} \|\mathcal{A}^*u^{k+1} + \mathcal{B}^*v^{k} - c\|= 0, \\[8pt]
\displaystyle\lim_{k\to\infty} (1-\tau)\|\mathcal{A}^*u^{k+1} + \mathcal{B}^*v^{k+1} - c\|= 0, \\[8pt]
\displaystyle \lim_{k\to\infty} \|w^{k+1} - w^k\|_{\mathcal{Q}+\text{Diag}\,(\mathcal{S}, \mathcal{T})} = 0.\\[8pt]
\end{array}
\end{equation}
For $\tau\in(0,1)$, since $\displaystyle \lim_{k\to\infty} \|\mathcal{A}^*u^{k+1} + \mathcal{B}^*v^{k+1} - c\|=0$, by using (\ref{lim_1}) we see that
\begin{equation*}\begin{array}{ll}
\displaystyle \lim_{k\to\infty} \|\mathcal{A}^*(u^{k+1} - u^k)\|\leq
\displaystyle \lim_{k\to\infty} (\|\mathcal{A}^*u^{k+1} + \mathcal{B}^*v^{k} - c\| + \|\mathcal{A}^*u^{k} + \mathcal{B}^*v^{k} - c\| ) = 0,\\[8pt]
\displaystyle \lim_{k\to\infty} \|\mathcal{B}^*(v^{k+1} - v^k)\|\leq
\displaystyle \lim_{k\to\infty}( \|\mathcal{A}^*u^{k+1} + \mathcal{B}^*v^{k+1} - c\| + \|\mathcal{A}^*u^{k+1} + \mathcal{B}^*v^{k} - c\|)  = 0,
\end{array}
\end{equation*}
 which implies $\displaystyle \lim_{k\to\infty} \|w^{k+1} - w^k\|_{\mathcal{Q}+\text{Diag}\,(\mathcal{S}+\sigma\mathcal{A}\mathcal{A}^*, \mathcal{T} + \sigma\mathcal{B}\mathcal{B}^*)} = 0$.
Therefore, for $\tau \in (0,1]$, we know that $\displaystyle \lim_{k\to\infty} \|w^{k+1} - w^k\|_{\mathcal{Q}+\text{Diag}\,(\mathcal{S}+(1-\tau)\sigma\mathcal{A}\mathcal{A}^*, \mathcal{T} + (1-\tau)\sigma\mathcal{B}\mathcal{B}^*)} = 0$. By condition (\ref{case1_assumption}) we can see that this implies either $\displaystyle \lim_{k\to\infty}\|u^{k+1} - u^k\| = 0$ or $\displaystyle \lim_{k\to\infty}\|v^{k+1} - v^k\| = 0$.
 Without loss of generality we assume that
$\displaystyle \lim_{k\to\infty} \|v^{k+1} - v^k\| = 0$. Thus,
\begin{equation}\label{Fineq}
\begin{array}{ll}
\displaystyle \lim_{k\to\infty} \|\mathcal{A}^*(u^{k+1} - u^k)\|& \leq
\displaystyle \lim_{k\to\infty}( \|\mathcal{A}^*u^{k+1} + \mathcal{B}^*v^{k} - c\| + \|\mathcal{A}^*u^{k} + \mathcal{B}^*v^{k-1} - c\| + \|\mathcal{B}^*(v^k-v^{k-1})\| )\\[8pt]
& = 0, \\[8pt]
\displaystyle \lim_{k\to\infty} \|u^{k+1} - u^k\|_{\mathcal{Q}_{11}}& \leq
\displaystyle \lim_{k\to\infty}( \|w^{k+1} - w^k\|_{\mathcal{Q}} + \|v^{k+1} - v^k\|_{\mathcal{Q}_{22}} ) = 0.
\end{array}
\end{equation}
Therefore, $\displaystyle \lim_{k\to\infty} \|u^{k+1} - u^k\|_{\mathcal{Q}_{11}+\mathcal{S}+\sigma\mathcal{A}\mathcal{A}^*} = 0$. This implies $\displaystyle \lim_{k\to\infty} \|u^{k+1} - u^k\|=0$ by the positive definiteness of $\mathcal{Q}_{11}+\mathcal{S}+\sigma\mathcal{A}\mathcal{A}^*$.

Now taking limits on both sides of (\ref{opt_1}) along the subsequence $\{(u^{k_i},v^{k_i},x^{k_i})\}$, and by using the closedness of the graphs of $\partial p$, $\partial q$ and the continuity of $\nabla \phi$, we obtain
\begin{equation*}\left\{\begin{array}{ll}
0\in F(u^\infty, v^\infty,x^\infty),\\[8pt]
\mathcal{A}^*u^\infty + \mathcal{B}^*v^\infty = c.
\end{array}\right.
\end{equation*}
This indicates that $(u^\infty,v^\infty)$ is an optimal solution to (\ref{opt}) and $x^\infty$ is the corresponding optimal multiplier.
Since $(u^\infty, v^\infty, x^{\infty})$ satisfies (\ref{kkt}), all the above arguments involving $(\bar{u},\bar{v},\bar{x})$ can be replaced by
$(u^\infty, v^\infty, x^{\infty})$. Thus  the subsequence $\{\Phi_{k_i}(u^\infty, v^\infty, x^{\infty}) \}$ converges to $0$ as $k_i\to \infty$.  Since $\{\Phi_{k_i}(u^\infty, v^\infty, x^{\infty}) \}$ is non-increasing, we obtain that
\begin{equation}\label{lim_phi}\begin{array}{ll}
\displaystyle\lim_{k\to\infty} \Phi_{k+1}(u^\infty, v^\infty, x^{\infty})= &\displaystyle\lim_{k\to\infty}~ (\tau\sigma)^{-1}\|x^{k+1} - x^\infty\|^2 + \|v^{k+1} - v^{\infty}\|^2_{\sigma\mathcal{B}\mathcal{B}^* + \mathcal{T}+\mathcal{Q}_{22}}+ \|u^{k+1} - u^{\infty}\|^2_{\mathcal{S}}\\[8pt]
& \quad\quad\quad +\|w^{k+1} - w^{\infty}\|_{\mathcal{Q}}^2 = 0.
\end{array}
\end{equation}
From this we can immediately get $\displaystyle\lim_{k\to\infty} x^{k+1} = x^\infty$ and $\displaystyle\lim_{k\to\infty} v^{k+1} = v^\infty$. Similar to inequality (\ref{Fineq}) we have that
$\displaystyle \lim_{k\to\infty} \sigma\|\mathcal{A}^*(u^{k+1} - u^{\infty})\| = 0$ and $\displaystyle \lim_{k\to\infty} \|u^{k+1} - u^{\infty}\|_{\mathcal{Q}_{11}} = 0$, which, together with (\ref{lim_phi}), imply that $\displaystyle \lim_{k\to\infty} \|u^{k+1} - u^{\infty}\| = 0$ by the positive definiteness of $\mathcal{Q}_{11}+\mathcal{S}+\sigma\mathcal{A}\mathcal{A}^*$.
Therefore,  the whole sequence $\{(u^k,v^k,x^k)\}$ converges to $(u^{\infty},v^\infty, x^\infty)$, the unique limit of the sequence. This completes the proof for the first case.
\bigskip

\noindent
{\bf (ii)} From the inequality (\ref{decrease_22}) and the optimality condition (\ref{VI})
we know that for any $k\geq 1$,
\begin{equation}\label{decrease_2}
\begin{array}{ll}
&
(\Psi_{k+1}(\bar{u}, \bar{v}, \bar{x}) + \Xi_{k+1}) - (\Psi_k(\bar{u}, \bar{v}, \bar{x}) +\Xi_k)\\[8pt]
\leq &
 -(\Gamma_{k+1}  +\min(1,1+\tau^{-1} - \tau)\sigma\|\mathcal{A}^*u^{k+1}   + \mathcal{B}^*v^{k+1} - c\|^2).
\end{array}
\end{equation}
By the assumptions  $\tau\in (0, \frac{1+\sqrt{5}}{2})$ and $\mathcal{M}\succeq 0$, we can obtain that $\Gamma_{k+1}\geq 0$ and $\min(1,1+\tau^{-1} - \tau)\geq 0$. Then both  $\{\Psi_{k+1}(\bar{u}, \bar{v}, \bar{x})\}$ and $\{\Xi_{k+1}\}$ are bounded. Thus, by a similar approach to case (i), we see that the sequence $\{(u^k,v^k,x^k)\}$ is bounded. Therefore, there exists a subsequence $\{(u^{k_i},v^{k_i}, x^{k_i})\}$ that converges to a cluster point, say $(u^{\infty},v^{\infty},x^{\infty})$. Next we will prove that
$(u^{\infty},v^{\infty})$ is optimal to (\ref{opt}) and $x^\infty$ is the corresponding optimal multiplier. The inequality (\ref{decrease_2}) also implies that
\begin{equation*}\begin{array}{ll}
\displaystyle \lim_{k\to\infty} \|x^{k+1} - x^k\| =  \displaystyle\lim_{k\to\infty}(\tau\sigma)^{-1}\|\mathcal{A}^*u^{k+1} + \mathcal{B}^*v^{k+1} - c\|= 0, \\[8pt]
\displaystyle \lim_{k\to\infty} \|w^{k+1} - w^k\|_{\mathcal{M}} = 0,\quad
\displaystyle \lim_{k\to\infty}\|\mathcal{B}(v^{k+1} - v^k)\| = 0.
\end{array}
\end{equation*}
By the relationship
\begin{equation*}\begin{array}{ll}
\displaystyle\lim_{k\to\infty} \|\mathcal{A}^*(u^{k+1}-u^k)\| & \leq \displaystyle\lim_{k\to\infty} (\|\mathcal{A}^*u^{k+1} + \mathcal{B}^*v^{k+1}-c\| + \|\mathcal{A}^*u^k + \mathcal{B}^*v^k - c\| + \|\mathcal{B}^*(v^{k+1} - v^k)\|)\\[8pt]
& = 0,
\end{array}
\end{equation*}
we can further get $\displaystyle \lim_{k\to\infty} \|w^{k+1} - w^k\|_{\mathcal{M} + \text{Diag}\,(\sigma\mathcal{A}\mathcal{A}^*, \sigma\mathcal{B}\mathcal{B}^*)}= 0$. Thus, by the condition (\ref{case2_assumption}), we can get that either $\displaystyle\lim_{k\to\infty} \|u^{k+1} - u^k\| = 0$
 or $\displaystyle\lim_{k\to\infty} \|v^{k+1} - v^k\| = 0$. Again similar to  case (i), we can see that in fact both of them would hold by the positive definiteness of $\frac{1}{4} \mathcal{Q}_{11}+\mathcal{S} + \sigma\mathcal{A}\mathcal{A}^* - \eta\mathcal{D}_1$ and $\frac{1}{4} \mathcal{Q}_{22}+\mathcal{T} + \sigma\mathcal{B}\mathcal{B}^* - \eta\mathcal{D}_{22}$. The remaining proof about the  convergence of the whole sequence $\{(u^k,v^k,x^k)\}$ follows exactly the same as in case (i).
This completes the proof for the second case.

\qed
\begin{remark}
In Theorem \ref{global_convergence},  for  $\tau\in(0,1]$, a sufficient condition for the convergence is
 \[\mathcal{Q}+\text{Diag}\,(\mathcal{S}+(1-\tau)\sigma\mathcal{A}\mathcal{A}^*,
\mathcal{T} + (1-\tau)\sigma\mathcal{B}\mathcal{B}^*)\succ 0,\]
and for $\tau\in [1, \frac{1+\sqrt{5}}{2})$, a sufficient condition for the convergence is
\begin{equation*}\begin{array}{ll}
\disp\frac{1}{4}\mathcal{Q} + \text{Diag}\,(\mathcal{S}-\eta\mathcal{D}_1, ~\mathcal{T} -\eta\mathcal{D}_2 )\succeq 0,\quad
\frac{1}{4}\mathcal{Q}+\text{Diag}\,(\mathcal{S}+\sigma\mathcal{A}\mathcal{A}^*-\eta\mathcal{D}_1, \mathcal{T} +\sigma\mathcal{B}\mathcal{B}^* -\eta\mathcal{D}_2)\succ 0.
\end{array}
\end{equation*}
\end{remark}

\begin{remark}\label{Q_plus_separable}
An interesting application of Theorem \ref{global_convergence} is for the linearly constrained convex optimization problem with a quadratically  coupled objective function of the form \[\phi(w) = \frac{1}{2}\langle w, \widetilde{\mathcal{Q}}w\rangle + f(u) + g(v),\] where $\widetilde{\mathcal{Q}}:\mathcal{U}\times\mathcal{V}\to\mathcal{U}\times\mathcal{V}$ is a self-adjoint positive semidefinite linear operator, $f: \mathcal{U}\to (-\infty, \infty)$ and $g: \mathcal{V}\to (-\infty, \infty)$ are two convex smooth functions with Lipschitz continuous gradients.
In this case, there exist four self-adjoint positive semidefinite operators $\Sigma_f, \widehat{\Sigma}_f: \mathcal{U}\to\mathcal{U}$ and $\Sigma_g, \widehat{\Sigma}_g:\mathcal{V}\to\mathcal{V}$
such that
\begin{equation*}
\Sigma_f\preceq \xi\preceq \widehat{\Sigma}_f \quad\forall \xi\in\partial^2 f(u), \,u\in\mathcal{U}
\end{equation*}
and
\begin{equation*}
\Sigma_g\preceq \zeta\preceq \widehat{\Sigma}_g \quad\forall \zeta\in\partial^2 g(v), \, v\in\mathcal{V},
\end{equation*}
where $\partial^2 f$ and $\partial^2 g$ are defined in (\ref{gen_Hessian}).
Then by letting $\mathcal{Q} = \widetilde{\mathcal{Q}} + \text{Diag}\,(\Sigma_f, \Sigma_g)$ in (\ref{phi_convex}) and $\mathcal{Q} + \mathcal{H} = \widetilde{\mathcal{Q}} + \text{Diag}\,(\widehat{\Sigma}_f, \widehat{\Sigma}_g)$ in (\ref{phi_majorize}),  we have $\eta = 0$ in (\ref{cross_control}). This implies that $\mathcal{M}\succeq 0$ always holds in (\ref{convergence_2}). Therefore,  for $\tau\in(0,\frac{1+\sqrt{5}}{2})$,
the conditions for the convergence can be equivalently written as
\begin{equation}\label{QP_assumption_1}\begin{array}{ll}
\widetilde{\mathcal{Q}}_{11} +\Sigma_f+ \mathcal{S} + \sigma\mathcal{A}\mathcal{A}^*\succ 0, ~\widetilde{\mathcal{Q}}_{22}+\Sigma_g + \mathcal{T} + \sigma\mathcal{B}\mathcal{B}^*\succ 0
\end{array}
\end{equation}
and
 \begin{equation}\label{QP_assumption_2}
\langle w, [\widetilde{\mathcal{Q}}+\text{Diag}\,(\Sigma_f + \mathcal{S}+\sigma\mathcal{A}\mathcal{A}^*, \Sigma_g + \mathcal{T} + \sigma\mathcal{B}\mathcal{B}^*)]w\rangle = 0 \Rightarrow \|u\|\|v\| = 0.
\end{equation}
A sufficient condition for ensuring (\ref{QP_assumption_1}) and (\ref{QP_assumption_2}) to hold is
\begin{equation}
\widetilde{\mathcal{Q}} + \text{Diag}\,(\Sigma_f + \mathcal{S} + \sigma\mathcal{A}\mathcal{A}^*, \Sigma_g + \mathcal{T} + \sigma\mathcal{B}\mathcal{B}^*)\succ 0.
\end{equation}
If $\widetilde{\mathcal{Q}} = 0$, i.e., if the objective function of the original problem (\ref{opt}) is separable, we will recover the convergence conditions given in \cite{LiSunToh_indefinite2014} for a majorized ADMM with semi-proximal terms.
\end{remark}

\subsection{The non-ergodic iteration complexity for general coupled objective functions }
In this section, we will present the non-ergodic iteration complexity for the majorized ADMM in terms of the KKT optimality condition.
\begin{theorem}\label{nonerg_complexity}
Suppose that the solution set of (\ref{opt}) is nonempty and Assumption \ref{cq} holds. Assume that
one of the following conditions holds:\\[8pt]
(i) $\tau\in(0,1]$ and
$ \mathcal{O}_1: = \disp\frac{1}{4}\mathcal{Q} + \text{Diag}\,(\mathcal{S}+(1-\tau)\sigma\mathcal{A}\mathcal{A}^*, \mathcal{T} + (1-\tau)\sigma\mathcal{B}\mathcal{B}^*)\succ 0$;\\[8pt]
(ii) $\tau\in (0,\frac{1+\sqrt{5}}{2})$, $\disp\frac{1}{4}\mathcal{Q} + \text{Diag}\,(\mathcal{S}- \eta\mathcal{D}_1, \mathcal{T}  - \eta\mathcal{D}_2)\succeq 0$ and
$\mathcal{O}_2: = \disp\frac{1}{4}\mathcal{Q} + \text{Diag}\,(\mathcal{S}+\sigma\mathcal{A}\mathcal{A}^* - \eta\mathcal{D}_1, \mathcal{T} + \sigma\mathcal{B}\mathcal{B}^* - \eta\mathcal{D}_2)\succ 0.$\\[8pt]
Then there exists a constant $C$ only depending on the initial point and the optimal solution set such that
the sequence $\{(u^k, v^k,x^k)\}$ generated by the majorized ADMM satisfies that for $k\geq 1$,
\begin{equation}\label{nonerg_main}
\displaystyle\min_{1\leq i\leq k} \{\text{dist}^2(0, F(u^{i+1}, v^{i+1}, x^{i+1}))
 + \|\mathcal{A}^*u^{i+1} + \mathcal{B}^*v^{i+1}-c\|^2\} \leq C/k.
 \end{equation}
 {
Furthermore,  for the limiting case we have that
 \begin{equation}\label{nonerg_limit}
 \displaystyle\lim_{k\to\infty} k(\min_{1\leq i\leq k} \{\text{dist}^2(0, F(u^{i+1}, v^{i+1}, x^{i+1}))
 + \|\mathcal{A}^*u^{i+1} + \mathcal{B}^*v^{i+1}-c\|^2 \})= 0.
\end{equation}
}
\end{theorem}
{\bf Proof.}
From the optimality condition for $(u^{k+1}, v^{k+1})$,  we know that
\begin{equation*}\begin{array}{ll}
&\left(\begin{array}{cc}
-(1-\tau)\sigma\mathcal{A}(\mathcal{A}^*u^{k+1} +\mathcal{B}^*v^{k+1} - c) - \sigma\mathcal{A}\mathcal{B}^*(v^{k} - v^{k+1})  - \mathcal{S}(u^{k+1} - u^k) + \mathcal{Q}_{12}(v^{k+1}-  v^k)\\
-(1-\tau)\sigma\mathcal{B}(\mathcal{A}^*u^{k+1} +\mathcal{B}^*v^{k+1} - c) -  \mathcal{T}(v^{k+1} - v^k)
\end{array}\right) \\[10pt]
&
- (\mathcal{Q}+\mathcal{H})(w^{k+1} - w^k)+ \nabla \phi(w^{k+1}) - \nabla \phi(w^k)\\[8pt]
 \in & F(u^{k+1}, v^{k+1}, x^{k+1}).
\end{array}
\end{equation*}
Therefore, we can obtain that
\begin{equation}\label{nonerg_1}\begin{array}{ll}
&\text{dist}^2(0, F(u^{k+1}, v^{k+1}, x^{k+1}))
 + \|\mathcal{A}^*u^{k+1} + \mathcal{B}^*v^{k+1}-c\|^2\\[8pt]

\leq & 5\|\sigma\mathcal{A}\mathcal{B}^*(v^{k+1} - v^k)\|^2 + 5(1-\tau)^2\sigma^2(\|\mathcal{A}\|^2+\|\mathcal{B}\|^2)\|\mathcal{A}^*u^{k+1} + \mathcal{B}^*v^{k+1} - c\|^2\\[8pt]
&
 + 5\| (\mathcal{Q} + \mathcal{H})(w^{k+1} - w^k)-\nabla\phi(w^{k+1}) +\nabla \phi(w^k)\|^2 + 5\|\mathcal{Q}_{12}(v^{k+1} - v^k)\|^2
 + 5\|\mathcal{T}(v^{k+1} - v^k)\|^2 \\[8pt]
 &
+ 5\|\mathcal{S}(u^{k+1} - u^k)\|^2+\|\mathcal{A}^*u^{k+1} + \mathcal{B}^*v^{k+1} - c\|^2\\[8pt]

\leq & 5\sigma\|\mathcal{A}\|^2\|v^{k+1} - v^k\|_{\sigma\mathcal{B}\mathcal{B}^*}^2
 + (5(1-\tau)^2\sigma^2(\|\mathcal{A}\|^2 + \|\mathcal{B}\|^2)+1)\|\mathcal{A}^*u^{k+1} + \mathcal{B}^*v^{k+1}-c\|^2 \\[8pt]
 & + 5\|\sqrt{\mathcal{Q}_{12}^*\mathcal{Q}_{12}}\|\|v^{k+1} - v^k\|_{\sqrt{\mathcal{Q}_{12}^*\mathcal{Q}_{12}}}^2
 + 5\|\mathcal{H}\|\|w^{k+1} - w^k\|_{\mathcal{H}}^2 + 5\|\mathcal{S}\|\|u^{k+1} - u^k\|^2_{\mathcal{S}}\\[8pt]
 &
 + 5\|\mathcal{T}\|\|v^{k+1} - v^k\|^2_{\mathcal{T}}\\[8pt]
\leq  & C_1\|w^{k+1} - w^k\|^2_{\widehat{\mathcal{O}}} + C_2\|\mathcal{A}^*u^{k+1} + \mathcal{B}^*v^{k+1}-c\|^2,

\end{array}
\end{equation}
where $$C_1 = 5\max(\sigma\|\mathcal{A}\|^2,\|\sqrt{\mathcal{Q}_{12}^*\mathcal{Q}_{12}}\|, \|\mathcal{H}\|,\|\mathcal{S}\|,\|\mathcal{T}\|), \quad C_2 = 5(1-\tau)^2\sigma^2(\|\mathcal{A}\|^2 + \|\mathcal{B}\|^2)+1,$$ 
$$\widehat{\mathcal{O}} = \mathcal{H} + \text{Diag}\,(\mathcal{S}, \mathcal{T} + \sigma\mathcal{B}\mathcal{B}^* + \sqrt{\mathcal{Q}_{12}^*\mathcal{Q}_{12}})$$ and the second inequality comes from the fact that there exists some $\mathcal{W}^k\in\text{conv}\{\partial^2\phi([w^{k-1}, w^k])\}$ such that
\begin{equation*}\begin{array}{ll}
 & \|(\mathcal{Q} + \mathcal{H})(w^{k+1} - w^k)-\nabla \phi(w^{k+1}) + \nabla \phi(w^k)\|^2\\[8pt]
= & \|(\mathcal{Q} + \mathcal{H}-\mathcal{W}^k)(w^{k+1} - w^k)\|^2
\leq  \|\mathcal{H}\|\|w^{k+1} - w^k\|^2_{\mathcal{H}}.
\end{array}
\end{equation*}
Next we will estimate the upper bounds for $\|w^{k+1} - w^k\|^2_{\widehat{\mathcal{O}}}$ and $\|\mathcal{A}^*u^{k+1} + \mathcal{B}^*v^{k+1} - c\|^2$ by only involving the initial point and the optimal solution set under the two different conditions.

First, assume condition (i) holds. For $\tau\in(0,1]$, by using (\ref{convergence_1}) in the proof of Theorem \ref{global_convergence},
we have that for $i\geq 1$,
\begin{equation*}\begin{array}{ll}
&\|w^{i+1} - w^i\|^2_{\frac{1}{4}\mathcal{Q} + \text{Diag}\,(\mathcal{S}, \mathcal{T})} + \sigma\|\mathcal{A}^*u^{i+1}  + \mathcal{B}^*v^i - c\|^2 + (1-\tau)\sigma\|\mathcal{A}^*u^{i+1}  + \mathcal{B}^*v^{i+1} -c\|^2 \\[8pt]
\leq & \Phi_i(\bar{u}, \bar{v}, \bar{x}) - \Phi_{i+1}(\bar{u}, \bar{v}, \bar{x}),
\end{array}
\end{equation*}
which, implies that,
\begin{equation*}\begin{array}{ll}
& \displaystyle\sum_{i=1}^k  (\|w^{i+1} - w^i\|^2_{\frac{1}{4}\mathcal{Q} + \text{Diag}\,(\mathcal{S}, \mathcal{T})} + \sigma\|\mathcal{A}^*u^{i+1}  + \mathcal{B}^*v^i - c\|^2 + (1-\tau)\sigma\|\mathcal{A}^*u^{i+1}  + \mathcal{B}^*v^{i+1} -c\|^2 \\[8pt]
\leq & \Phi_1(\bar{u}, \bar{v}, \bar{x}) - \Phi_{k+1}(\bar{u}, \bar{v}, \bar{x})\leq \Phi_1(\bar{u}, \bar{v}, \bar{x}).
\end{array}
\end{equation*}
This shows that
\begin{equation}\label{nonerg_7}\begin{array}{ll}
\displaystyle\sum_{i=1}^k  \|w^{i+1} - w^i\|^2_{\frac{1}{4}\mathcal{Q} + \text{Diag}(\mathcal{S}, \mathcal{T})} \leq \Phi_1(\bar{u}, \bar{v}, \bar{x}),\\[8pt]
\displaystyle\sum_{i=1}^k  \sigma\|\mathcal{A}^*u^{i+1}  + \mathcal{B}^*v^{i} -c\|^2 \leq \Phi_1(\bar{u}, \bar{v}, \bar{x}),\\[8pt]
\displaystyle\sum_{i=1}^k  (1-\tau)\sigma\|\mathcal{A}^*u^{i+1} + \mathcal{B}^*v^{i+1} - c\|^2 \leq \Phi_1(\bar{u}, \bar{v}, \bar{x}).
\end{array}
\end{equation}
From the above three inequalities we can also get that
\begin{equation*}\begin{array}{ll}
(1-\tau)\displaystyle\sum_{i=1}^k \|u^{i+1} - u^i\|^2_{\sigma\mathcal{A}\mathcal{A}^*}& \leq (1-\tau)\displaystyle\sum_{i=1}^k (2\sigma\|\mathcal{A}^*u^{i+1} + \mathcal{B}^*v^i - c\|^2 + 2\sigma\|\mathcal{A}^*u^{i} + \mathcal{B}^*v^{i} - c\|^2) \\[8pt]
&\leq 2(2-\tau)\Phi_1(\bar{u}, \bar{v}, \bar{x}),\\[8pt]
(1-\tau)\displaystyle\sum_{i=1}^k \|v^{i+1} - v^i\|^2_{\sigma\mathcal{B}\mathcal{B}^*}& \leq (1-\tau)\displaystyle\sum_{i=1}^k (2\sigma\|\mathcal{A}^*u^{i+1} + \mathcal{B}^*v^i - c\|^2 + 2\sigma\|\mathcal{A}^*u^{i+1} + \mathcal{B}^*v^{i+1} - c\|^2) \\[8pt]
&
\leq 2(2-\tau)\Phi_1(\bar{u}, \bar{v}, \bar{x}).
\end{array}
\end{equation*}
With the notation of operator $\mathcal{O}_1$ we have that
\begin{equation}\label{nonerg_3}\begin{array}{ll}
\displaystyle\sum_{i=1}^k  \|w^{i+1} - w^i\|^2_{\mathcal{O}_1} &=\displaystyle\sum_{i=1}^k  \|w^{i+1} - w^i\|^2_{\frac{1}{4}\mathcal{Q} + \text{Diag}(\mathcal{S}, \mathcal{T})} + \displaystyle\sum_{i=1}^k \|w^{i+1} - w^i\|^2_{(1-\tau)\text{Diag}(\sigma\mathcal{A}\mathcal{A}^*, \mathcal{B}\mathcal{B}^*)}\\[8pt]
& \leq (9-4\tau)\Phi_1(\bar{u}, \bar{v}, \bar{x}).
\end{array}
\end{equation}
If $\tau\in(0,1)$, we further have that
\begin{equation}\label{nonerg_8}\begin{array}{ll}
\displaystyle\sum_{i=1}^k  \|\mathcal{A}^*u^{i+1} + \mathcal{B}^*v^{i+1} - c\|^2 \leq (1-\tau)^{-1}\sigma^{-1}\Phi_1(\bar{u}, \bar{v}, \bar{x}).
\end{array}
\end{equation}
If $\tau = 1$, by the condition that $\mathcal{O}_1 = \disp\frac{1}{4}\mathcal{Q} + \text{Diag}\,(\mathcal{S}, \mathcal{T})\succ 0$, we have that
\begin{equation}\label{nonerg_9}\begin{array}{ll}
\displaystyle\sum_{i=1}^k \|\mathcal{A}^*u^{i+1} + \mathcal{B}^*v^{i+1} - c\|^2 & \leq \displaystyle\sum_{i=1}^k (2\|\mathcal{A}^*u^{i+1} + \mathcal{B}^*v^{i} - c\|^2 + 2\|v^{i+1} - v^i\|^2_{\mathcal{B}\mathcal{B}^*})\\[8pt]
& \leq \displaystyle\sum_{i=1}^k (2\|\mathcal{A}^*u^{i+1} + \mathcal{B}^*v^{i} - c\|^2 \\[8pt]
&\quad \quad\quad+ 2\|\mathcal{O}_1^{-\frac{1}{2}}\text{Diag}\,(0,\mathcal{\mathcal{B}\mathcal{B}^*})\mathcal{O}_1^{-\frac{1}{2}}\|\|w^{i+1} - w^i\|^2_{\mathcal{O}_1})\\[8pt]
&
\leq (2\sigma^{-1} + (18-8\tau)\|\mathcal{O}_1^{-\frac{1}{2}}\text{Diag}\,(0,\mathcal{\mathcal{B}\mathcal{B}^*})\mathcal{O}_1^{-\frac{1}{2}}\|)\Phi_1(\bar{u}, \bar{v}, \bar{x}),
\end{array}
\end{equation}
where the second inequality is obtained by the fact that for any $\xi$, a  self-adjoint positive definite operator $\mathcal{G}$ with square root $\mathcal{G}^{\frac{1}{2}}$ and a self-adjoint positive semidefinite operator $\widehat{G}$ defined in the same Hilbert space, it always holds that
\begin{equation*}\begin{array}{ll}
\|\xi\|^2_{\widehat{\mathcal{G}}} = \langle \xi, \widehat{\mathcal{G}}\xi\rangle  = \langle \xi, (\mathcal{G}^{\frac{1}{2}}\mathcal{G}^{-\frac{1}{2}})\widehat{\mathcal{G}}(\mathcal{G}^{-\frac{1}{2}}\mathcal{G}^{\frac{1}{2}})\xi\rangle
& = \langle \mathcal{G}^{\frac{1}{2}}\xi, (\mathcal{G}^{-\frac{1}{2}}\widehat{\mathcal{G}}\mathcal{G}^{-\frac{1}{2}})\mathcal{G}^{\frac{1}{2}}\xi\rangle \\[8pt]
&
\leq \|\mathcal{G}^{-\frac{1}{2}}\widehat{\mathcal{G}}\mathcal{G}^{-\frac{1}{2}}\|\|\xi\|^2_{\mathcal{G}}.
\end{array}
\end{equation*}
Therefore, by using (\ref{nonerg_1}), (\ref{nonerg_3}) and the positive definiteness of operator $\mathcal{O}_1$, we know that
\begin{equation*}\begin{array}{ll}
&\displaystyle\min_{1\leq i\leq k} \{\text{dist}^2(0, F(u^{i+1}, v^{i+1}, x^{i+1}))
 + \|\mathcal{A}^*u^{i+1} + \mathcal{B}^*v^{i+1}-c\|^2\}\\[8pt]
\leq & (\displaystyle\sum_{i=1}^k(\text{dist}^2(0, F(u^{i+1}, v^{i+1}, x^{i+1}))
 + \|\mathcal{A}^*u^{i+1} + \mathcal{B}^*v^{i+1}-c\|^2))/k
\leq  C\Phi_1(\bar{u}, \bar{v}, \bar{x})/k,
\end{array}
 \end{equation*}
where
\begin{equation*}C =\left\{\begin{array}{ll}
C_1(9-4\tau)\|\mathcal{O}_1^{-\frac{1}{2}}\widehat{\mathcal{O}}\mathcal{O}_1^{-\frac{1}{2}}\| + C_2(1-\tau)^{-1}\sigma^{-1}, \quad\tau\in(0,1),\\[8pt]
C_1(9-4\tau)\|\mathcal{O}_1^{-\frac{1}{2}}\widehat{\mathcal{O}}\mathcal{O}_1^{-\frac{1}{2}}\| + C_2(2\sigma^{-1} + (18-8\tau)\|\mathcal{O}_1^{-\frac{1}{2}}\text{Diag}\,(0,\mathcal{\mathcal{B}\mathcal{B}^*})\mathcal{O}_1^{-\frac{1}{2}}\|)
,\quad \tau = 1.
\end{array}\right.
\end{equation*}
To prove the limiting case (\ref{nonerg_limit}),  by using inequalities (\ref{nonerg_3}), (\ref{nonerg_8}), (\ref{nonerg_9}) and \cite[Lemma 2.1]{LiSunToh_indefinite2014}, we have that
\begin{equation*}
\min_{1\leq i\leq k}\|w^{i+1} - w^i\|^2_{\mathcal{O}_1} = o(1/k),\quad
\min_{1\leq i\leq k}\|\mathcal{A}^*u^{i+1} + \mathcal{B}^*v^{i+1} - c\|^2 = o(1/k),
\end{equation*}
which, together with (\ref{nonerg_1}), imply that
\begin{equation*}\begin{array}{ll}
&\displaystyle\lim_{k\to\infty} k (\min_{1\leq i\leq k} \{\text{dist}^2(0, F(u^{i+1}, v^{i+1}, x^{i+1}))
 + \|\mathcal{A}^*u^{i+1} + \mathcal{B}^*v^{i+1}-c\|^2\})\\[8pt]
\leq &\displaystyle\lim_{k\to\infty} k (\min_{1\leq i\leq k}  \{C_1\|\mathcal{O}_1^{-\frac{1}{2}}\widehat{\mathcal{O}}\mathcal{O}_1^{-\frac{1}{2}}\|\|w^{i+1} - w^i\|^2_{\mathcal{O}_1} + C_2\|\mathcal{A}^*u^{i+1} + \mathcal{B}^*v^{i+1} - c\|^2\}) = 0.
 \end{array}
\end{equation*}
Now, we complete the proof of the conclusions under condition (i). Next, we consider the case under condition (ii). The proof of this part is similar to the case under condition (i). For  $\tau\in(0,\frac{1+\sqrt{5}}{2})$,
let $\rho(\tau) = \min(\tau,1+\tau-\tau^2)$. 
We know
from (\ref{convergence_2}) that for $\tau\in(0,\frac{1+\sqrt{5}}{2})$ and any $k\geq 1$,
\begin{eqnarray*}
&&\hspace{-0.7cm}
\displaystyle\sum_{i=1}^k  \|w^{i+1} - w^i\|^2_{\frac{1}{4}\mathcal{Q} + \text{diag}\,(\mathcal{S} - \eta\mathcal{D}_1,\mathcal{T}  - \eta\mathcal{D}_2)} 
+ \rho(\tau) \|v^{i+1} - v^i\|^2_{\sigma
\mathcal{B}\mathcal{B}^*}+  \frac{\rho(\tau)}{\tau}\sigma\|\mathcal{A}^*u^{i+1} + \mathcal{B}^*v^{i+1} - c\|^2 
\\[5pt]
&\leq &  (\Psi_1(\bar u, \bar v, \bar x) +\Xi_1) - (\Psi_{k+1}(\bar u, \bar v, \bar x) + \Xi_{k+1})\leq \Psi_1(\bar u, \bar v, \bar x) + \Xi_1.
\end{eqnarray*}
Thus by the positive semidefiniteness of $\frac{1}{4}\mathcal{Q} + \text{Diag}\,(\mathcal{S}- \eta\mathcal{D}_1, \mathcal{T} - \eta\mathcal{D}_2)$, we can get that
\begin{equation}\label{nonerg_4}\begin{array}{ll}
\displaystyle\sum_{i=1}^k  \|w^{i+1} - w^i\|^2_{\frac{1}{4}\mathcal{Q} + \text{diag}\,(\mathcal{S} - \eta\mathcal{D}_1,\mathcal{T}  - \eta\mathcal{D}_2)} \leq \Psi_1(\bar u, \bar v, \bar x) + \Xi_1,\\[8pt]
\displaystyle\sum_{i=1}^k  \|v^{i+1} - v^i\|^2_{\sigma\mathcal{B}\mathcal{B}^*} \leq (\Psi_1(\bar u, \bar v, \bar x)+\Xi_1)/\rho(\tau),\\[8pt]
\displaystyle\sum_{i=1}^k  \sigma\|\mathcal{A}^*u^{i+1} + \mathcal{B}^*v^{i+1} - c\|^2 \leq \tau (\Psi_1(\bar u, \bar v, \bar x)+\Xi_1)/\rho(\tau),
\end{array}
\end{equation}
which, implies that,
\begin{equation}\label{nonerg_5}\begin{array}{ll}
\displaystyle\sum_{i=1}^k  \|u^{i+1} - u^i\|^2_{\sigma\mathcal{A}\mathcal{A}^*} & \leq \displaystyle\sum_{i=1}^k(3\sigma\|\mathcal{A}^*u^{i+1} + \mathcal{B}^*v^{i+1} - c\|^2 + 3\sigma\|\mathcal{A}^*u^i + \mathcal{B}^*v^i - c\|^2 + 3\|v^{i+1} - v^i\|_{\sigma\mathcal{B}\mathcal{B}^*}^2)\\[8pt]
&
\leq (6\tau+2)(\Psi_1(\bar u, \bar v, \bar x)+\Xi_1)/\rho(\tau).
\end{array}
\end{equation}
Combining (\ref{nonerg_4}) and (\ref{nonerg_5}) one can find that
\begin{equation}\label{nonerg_6}\begin{array}{ll}
\displaystyle\sum_{i=1}^k  \|w^{i+1} - w^i\|^2_{\mathcal{O}_2}&  = \displaystyle\sum_{i=1}^k  \|w^{i+1} - w^i\|^2_{\frac{1}{4}\mathcal{Q} + \text{diag}\,(\mathcal{S} - \eta\mathcal{D}_1,\mathcal{T}  - \eta\mathcal{D}_2)}  + \displaystyle\sum_{i=1}^k  \|w^{i+1} - w^i\|^2_{\text{Diag}\,(\sigma\mathcal{A}\mathcal{A}^*, \sigma\mathcal{B}\mathcal{B}^*)}\\[8pt]
& \leq (1+ (6\tau+3)/\rho(\tau))(\Psi_1(\bar u, \bar v, \bar x) +\Xi_1).
\end{array}
\end{equation}
Therefore, by using (\ref{nonerg_1}), (\ref{nonerg_4}), (\ref{nonerg_6}), and recalling the positive definiteness of operator $\mathcal{O}_2$, we finally have that
\begin{equation*}\begin{array}{ll}
& \displaystyle\min_{1\leq i\leq k} \{\text{dist}^2(0,F(u^{i+1}, v^{i+1}, x^{i+1}))
 + \|\mathcal{A}^*u^{i+1} + \mathcal{B}^*v^{i+1}-c\|^2\}\\[8pt]
\leq & \big(\displaystyle\sum_{i=1}^k (\text{dist}^2(0,F(u^{i+1}, v^{i+1}, x^{i+1}))
 + \|\mathcal{A}^*u^{i+1} + \mathcal{B}^*v^{i+1}-c\|^2\big)/k\\[8pt]
\leq & C'(\Psi_1(\bar u, \bar v, \bar x)+\Xi_1)/k,
\end{array}
\end{equation*}
where $C' = C_1\|\mathcal{O}_2^{-\frac{1}{2}}\widehat{\mathcal{O}}\mathcal{O}_2^{-\frac{1}{2}}\|(1+ (6\tau+3)/\rho(\tau))+ C_2\sigma^{-1}\tau/\rho(\tau)$.
The limiting property (\ref{nonerg_limit}) can be derived in the same way as for the case under condition (i).
This completes the proof of Theorem \ref{nonerg_complexity}.
\qed
\begin{remark}\label{nonerg_remark}
Theorem \ref{nonerg_complexity} gives the non-ergodic complexity of the
KKT optimality condition, which 
does not seem to be known even for the classic ADMM with separable objective functions. For the latter, Davis and Yin~\cite{Yin2014}  provided related  non-ergodic iteration complexity results for the primal feasibility and the objective functions and constructed an interesting example to show
for any $\alpha>1/2$,  there exists an initial point such that the sequence $\{(u^k,v^k)\}$ generated by the classic ADMM $(\tau = 1)$
satisfies  $\|\mathcal{A}^*u^{k+1} + \mathcal{B}^*v^{k+1} - c\| \geq 1/k^{\alpha}$. 
This implies that the non-ergodic iteration complexity results presented in 
Theorem \ref{nonerg_complexity} 
 in terms of the  KKT optimality condition may be optimal.
 \end{remark}

\subsection{The ergodic iteration complexity for general coupled objective functions}
In this section, we will discuss the ergodic iteration complexity of  the majorized ADMM for solving problem (\ref{opt}). For $k = 1,2,\cdots,$ denote
 \begin{equation*}\begin{array}{ll}
\hat{x}^k = \displaystyle\frac{1}{k}\displaystyle\sum_{i=1}^k \tilde{x}^{i+1}, \quad
\hat{u}^k = \frac{1}{k}\displaystyle\sum_{i=1}^k u^{i+1}, \quad \hat{v}^k = \frac{1}{k}\displaystyle\sum_{i=1}^k v^{i+1}, \quad \hat{w}^k = (\hat{u}^k, \hat{v}^k)
\end{array}
\end{equation*}
and
\begin{equation*}\left\{\begin{array}{ll}
\Lambda_{k+1} = \|u^{k+1}_e\|^2_{\mathcal{D}_1 + \mathcal{S}} +\|v^{k+1}_e\|^2_{\mathcal{D}_2 +\mathcal{T}+\mathcal{Q}_{22}+\sigma\mathcal{B}\mathcal{B}^*}  +(\tau\sigma)^{-1}\|x^{k+1}\|^2, \\[8pt]
\overline\Lambda_{k+1} = \Lambda_{k+1} + \Xi_{k+1} + \|w^{k+1}_e\|^2_{\mathcal{Q}} + \max(1-\tau, 1-\tau^{-1})\sigma\|\mathcal{A}^*u^{k+1} + \mathcal{B}^*v^{k+1} - c\|^2.

\end{array}\right.
\end{equation*}

\begin{theorem}\label{erg_main}
Suppose that 
$\mathcal{S}$ and $\mathcal{T}$ are chosen such that \[\mathcal{Q}_{11} + \sigma\mathcal{A}\mathcal{A}^*+\mathcal{S}\succ 0, \quad \mathcal{Q}_{22} + \sigma\mathcal{B}\mathcal{B}^*+\mathcal{T}\succ 0.\]
Assume that either (a) $\tau\in(0,1]$ and (\ref{case1_assumption}) holds    or (b) $\tau\in(0,\frac{1+\sqrt{5}}{2})$ and (\ref{convergence_2}) and (\ref{case2_assumption})  hold.
Then there exist constants $D_1$, $D_2$ and $D_3$ that only depending on the initial point and the optimal solution set such that for $k\geq 1$, the following conclusions hold:

\noindent
(i) \begin{equation}\label{erg_pri}
 \|\mathcal{A}^*\hat{u}^{k} + \mathcal{B}^*\hat{v}^{k} - c\|\leq D_1/k.
\end{equation}

\noindent
(ii) For any $(u,v,x)\in {B}_k:  = \{(u,v,x)\in \mathcal{U}\times\mathcal{V}\times\mathcal{X}|\|(u,v,x) - (\hat{u}^k, \hat{v}^k, \hat{x}^k)\| \leq 1\}$,
\begin{equation}\label{complexity_decrease_1}\begin{array}{ll}
 &(p(\hat{u}_k) + q(\hat{v}_k) ) -( p(u) + q(v) )+\langle \hat{w}^{k}-w,\nabla\phi(w)\rangle+ \langle   \hat{u}^{k} - u, \mathcal{A}x\rangle  + \langle  \hat{v}^{k} - v , \mathcal{B}x\rangle\\[8pt]
&
 - \langle  \hat{x}^{k} - x , \mathcal{A}^*u + \mathcal{B}^*v - c\rangle \leq D_2/k.
 \end{array}
\end{equation}

\noindent
(iii) For case (b), if we further assume 
that  $\mathcal{S}- \eta\mathcal{D}_1\succeq 0$ and $\mathcal{T}- \eta\mathcal{D}_2\succeq 0$, then
\begin{equation}\label{erg_obj}
|\theta(\hat{u}^{k} ,\hat{v}^{k}) - \theta(\bar{u}, \bar{v}) | \leq  D_3/k.
\end{equation}
The inequality (\ref{erg_obj}) holds for case (a) without additional assumptions.
\end{theorem}
{\bf Proof.}
{\bf (i)} Under the conditions for case (a), the inequality (\ref{convergence_1}) indicates that $\{\Phi_{k+1}(\bar{u}, \bar{v}, \bar{x})\}$ is a non-increasing sequence, which implies that
\begin{equation*}
(\tau\sigma)^{-1}\|x^{k+1} - \bar{x}\|^2 \leq \Phi_{k+1}(\bar{u}, \bar{v}, \bar{x}) \leq \Phi_1(\bar{u}, \bar{v}, \bar{x}).
\end{equation*}
Similarly under the conditions for case (b),  we can get from (\ref{decrease_2})  that
\begin{equation*}
(\tau\sigma)^{-1}\|x^{k+1} - \bar{x}\|^2\leq \Psi_{k+1}(\bar{u}, \bar{v}, \bar{x})+\Xi_{k+1}\leq \Psi_1(\bar{u}, \bar{v}, \bar{x}) + \Xi_1.
\end{equation*}
Therefore, in terms of the ergodic primal feasibility, we have that
\begin{equation}\label{erg_1}\begin{array}{ll}
\|\mathcal{A}^*\hat{u}^{k} + \mathcal{B}^*\hat{v}^{k} - c\|^2 & = \|\frac{1}{k}\displaystyle\sum_{i=1}^k (\mathcal{A}^*{u}^{i+1} + \mathcal{B}^*{v}^{i+1} - c)\|^2 \\[8pt]
&
 = \|(\tau\sigma)^{-1}\displaystyle\sum_{i=1}^k (x^{i+1} - x^i)\|^2/k^2\\[8pt]
&
 = \|(\tau\sigma)^{-1}(x^{k+1} - x^1)\|^2/k^2\\[8pt]
 &
\leq 2\|(\tau\sigma)^{-1}(x^{k+1} - \bar x)\|^2/k^2 + 2\|(\tau\sigma)^{-1}(x^1 - \bar x)\|^2/k^2
\leq C_3/k^2,
\end{array}
\end{equation}
where for case (a), $C_3 = 2(\tau\sigma)^{-1}\Phi_1(\bar{u}, \bar{v}, \bar{x}) +2\|(\tau\sigma)^{-1}(x^1 - \bar x)\|^2$ and for case (b), $C_3 = 2(\tau\sigma)^{-1}(\Psi_1(\bar{u}, \bar{v}, \bar{x}) + \Xi_1) +2\|(\tau\sigma)^{-1}(x^1 - \bar x)\|^2$. Then by taking the square root  on inequality (\ref{erg_1}),
we can obtain (\ref{erg_pri}).
\bigskip

\noindent
{\bf (ii)} First, assume that the conditions for case (a) hold. Then the
inequality (\ref{decrease_11}) implies that for $i\geq 1$,
 \begin{equation*}
\begin{array}{ll}
&(p(u^{i+1}) + q(v^{i+1}) ) -(p(u) + q(v) ) +\langle w^{i+1}-w,\nabla\phi(w)\rangle+ \langle   u^{i+1} - u, \mathcal{A}x\rangle  + \langle  v^{i+1} - v , \mathcal{B}x\rangle\\[8pt]
&
 - \langle  \tilde{x}^{i+1} - x , \mathcal{A}^*u + \mathcal{B}^*v - c\rangle \\[8pt]

\leq  & -\frac{1}{2}(\Phi_{i+1}({u}, {v}, {x}) - \Phi_i({u}, {v}, {x})).
\end{array}
\end{equation*}
Thus,
 summing up the above inequalities over $i = 1, \cdots,k$ and by using the convexity of functions $p$ and $q$, we can obtain that
\begin{equation}\label{vi_complexity}\begin{array}{ll}
p(\hat{u}_k) + q(\hat{v}_k) -  (p(u) + q(v))  +\langle \hat{w}^{k}-w,\nabla\phi(w)\rangle+ \langle   \hat{u}^{k} - u, \mathcal{A}x\rangle  + \langle  \hat{v}^{k} - v , \mathcal{B}x\rangle\\[8pt]

- \langle  \hat{x}^k -x, \mathcal{A}^*u + \mathcal{B}^*v - c\rangle

 \leq (\Phi_{1}(u,v,x) - \Phi_{k+1}(u,v,x))/2k \leq \Phi_1(u,v,x)/2k.
\end{array}
\end{equation}
Next, we will provide an explicit bound of $\Phi_1(u,v,x)$ for $(u,v,x)\in {B}_k$
 that only depends on the initial point and the optimal solution set.
Since $\{\Phi_{i+1}(\bar u,\bar v,\bar x)\}$ is non-increasing with $i$,  for $i\geq 1$ we have
\begin{equation*}
 (\tau \sigma)^{-1} \|\bar x - x^{i+1}
\|^2 + \|\bar u - u^{i+1} \|^2_{\mathcal{D}_1 + \mathcal{S} }  +
\|\bar v - v^{i+1}
\|^2_{\mathcal{Q}_{22} + \mathcal{D}_{2} + \mathcal{T} + \sigma\mathcal{B}\mathcal{B}^*} + \frac{1}{2}\|\bar w  - w^{i+1}\|^2_{\mathcal{Q}}
\leq   \Phi_1(\bar{u}, \bar{v}, \bar{x}).
\end{equation*}
Thus, summing up the above inequalities from $i = 1$ to $k$ and by applying Jensen's inequality to the convex function $\|\cdot\|^2$, we can get that
 \begin{equation}\label{ergodic_2}
 (\tau \sigma)^{-1} \|\bar x - \frac{1}{k}\displaystyle\sum_{i=1}^k x^{i+1}
\|^2   + \|\bar u - \hat u^k \|^2_{\mathcal{D}_1 + \mathcal{S} }  +
\|\bar v - \hat v^k
\|^2_{\mathcal{Q}_{22} + \mathcal{D}_{2} + \mathcal{T} + \sigma\mathcal{B}\mathcal{B}^*} + \frac{1}{2}\|\bar w  - \hat w^k\|^2_{\mathcal{Q}}
\leq  \Phi_1(\bar{u}, \bar{v}, \bar{x}).
\end{equation}
Recall that $\tilde x^{i+1}  = x^i + \sigma(\mathcal{A}^*u^{i+1} +\mathcal{B}^*v^{i+1} - c) = x^i +\tau^{-1}(x^{i+1} - x^i)$. Then we also have
\begin{equation*}
\sum_{i=1}^k\tilde x^{i+1} = \sum_{i=1}^k x^i + \tau^{-1}(x^{k+1} - x^1) = \sum_{i=1}^k x^{i+1} + (\tau^{-1} -1)(x^{k+1} - x^1),
\end{equation*}
which, implies that for $k\geq 1$,
\begin{equation}\label{ergodic_3}\begin{array}{ll}
\|\bar{x} - \hat{x}^k\|^2 &=  \|\bar{x} - \disp\frac{1}{k}(\displaystyle\sum_{i=1}^k x^{i+1} + (\tau^{-1} -1)(x^{k+1} - x^1))\|^2\\[8pt]
&\leq 2\|\bar{x} - \disp\frac{1}{k}\displaystyle\sum_{i=1}^k x^{i+1}
\|^2  + 2(\tau^{-1} - 1)^2\|\frac{1}{k}(x^{k+1} - x^1)\|^2\\[8pt]
&\leq 2\|\bar{x} - \disp\frac{1}{k}\displaystyle\sum_{i=1}^k x^{i+1}
\|^2  + 4(\tau^{-1} - 1)^2\|x^{k+1} -\bar{x}\|^2 + 4(\tau^{-1} - 1)^2\|x^{1} -\bar{x}\|^2\\[8pt]
&\leq 2\|\bar{x} - \disp\frac{1}{k}\displaystyle\sum_{i=1}^k x^{i+1}
\|^2  + 8(\tau^{-1} - 1)^2\Phi_1(\bar{u},\bar{v}, \bar{x}).
\end{array}
\end{equation}
Therefore,  by (\ref{ergodic_2}) and (\ref{ergodic_3}), we can obtain that
\begin{equation}\label{ergodic_4}\begin{array}{ll}
&(\tau \sigma)^{-1} \|\bar x - \hat{x}^k
\|^2 + \|\bar u - \hat{u}^k \|^2_{\mathcal{D}_1 + \mathcal{S} }  + \|\bar v - \hat{v}^k\|^2_{\mathcal{Q}_{22} + \mathcal{D}_{2} + \mathcal{T} + \sigma\mathcal{B}\mathcal{B}^*}  + \disp\frac{1}{2}\| \bar w - \hat{w}^k\|^2_{\mathcal{Q}}\\[8pt]
\leq & 2(\tau\sigma)^{-1}\|\bar{x} - \disp\frac{1}{k}\displaystyle\sum_{i=1}^k x^{i+1}
\|^2   + \|\bar u - \hat{u}^k \|^2_{\mathcal{D}_1 + \mathcal{S} }  + \|\bar v - \hat{v}^k\|^2_{\mathcal{Q}_{22} + \mathcal{D}_{2} + \mathcal{T} + \sigma\mathcal{B}\mathcal{B}^*}  + \frac{1}{2}\| \bar w - \hat{w}^k\|^2_{\mathcal{Q}}\\[8pt]
& +8(\tau^{-1} - 1)^2\Phi_1(\bar{u},\bar{v}, \bar{x})\\[8pt]
\leq & (2 + 8(\tau^{-1}-1)^2)\Phi_1(\bar{u},\bar{v}, \bar{x}).
\end{array}
\end{equation}
In addition, we also know from (\ref{erg_1})  that for $k\geq 1$,
\begin{equation}\label{ergodic_5}\begin{array}{ll}
\|\mathcal{A}^*\hat{u}^{k} + \mathcal{B}^*\hat{v}^{k} -c\|^2 \leq C_3/k^2\leq C_3.
\end{array}
\end{equation}
Now putting inequalities (\ref{ergodic_4}) and (\ref{ergodic_5}) together,  we see that for any {$(u,v,x)\in {B}_k,$}
\begin{equation*}
\begin{array}{ll}
& \disp(\tau \sigma)^{-1} \|x^1 - x\|^2+ \|u^1 -
u\|^2_{\mathcal{D}_1 + \mathcal{S}} + \|v^1-v\|^2_{\mathcal{Q}_{22} + \mathcal{D}_{2} +\mathcal{T} }  + \frac{1}{2}\|w^{1} - w\|^2_{\mathcal{Q}}  + \sigma \|\mathcal{A}^*u + \mathcal{B}^*v^1 - c\|^2 \\[8pt]
= & (\tau \sigma)^{-1} \|(x^1 - \bar{x}) + (\bar{x} - \hat{x}^k) + (\hat{x}^k - x)\|^2+ \|(u^1 -\bar{u}) + (\bar{u} - \hat{u}^k) + (\hat{u}^k - u)\|^2_{\mathcal{D}_1 + \mathcal{S}} \\[8pt]
&
+ \disp\|(v^1-\bar v) + (\bar{v} - \hat{v}^k) + (\hat{v}^k - v)\|^2_{\mathcal{Q}_{22} + \mathcal{D}_{2} +\mathcal{T} }  + \frac{1}{2}\|(w^{1} - \bar{w}) +( \bar{w} - \hat{w}^k) + (\hat{w}^k - w)\|^2_{\mathcal{Q}} \\[8pt]
&
 \disp+ \sigma \|\mathcal{A}^*(u - \hat{u}^k)  + \mathcal{B}^*(v^1 -\bar{v}) + \mathcal{B}^*(\bar{v}- \hat{v}^k)+( \mathcal{A}^*\hat{u}^k + \mathcal{B}^*\hat{v}^k - c)\|^2\\[8pt]
\leq & 3[(\tau \sigma)^{-1} \|x^1 - \bar x
\|^2 + \|u^1 - \bar u \|^2_{\mathcal{D}_1 + \mathcal{S} }  + \|v^1 - \bar v\|^2_{\mathcal{Q}_{22} + \mathcal{D}_{2} + \mathcal{T} + 2\sigma\mathcal{B}\mathcal{B}^*}  + \frac{1}{2}\| w^1 - \bar w\|^2_{\mathcal{Q}}] \\[8pt]
&\disp + 3[(\tau \sigma)^{-1} \|\bar x - \hat{x}^k
\|^2  + \|\bar u - \hat{u}^k \|^2_{\mathcal{D}_1 + \mathcal{S} }  + \|\bar v - \hat{v}^k\|^2_{\mathcal{Q}_{22} + \mathcal{D}_{2} + \mathcal{T} + \sigma\mathcal{B}\mathcal{B}^*}  + \frac{1}{2}\| \bar w - \hat{w}^k\|^2_{\mathcal{Q}}]\\[8pt]
&    + 3[(\tau \sigma)^{-1} \|\hat{x}^k - x\|^2 + \|
\hat{u}^k - u\|^2_{\mathcal{D}_1 + \mathcal{S}+2\sigma \mathcal{A}\mathcal{A}^*}  + \|
\hat{v}^k - v\|^2_{\mathcal{Q}_{22} + \mathcal{D}_2 + \mathcal{T}} + \frac{1}{2}\|\hat{w}^k - w\|^2_{\mathcal{Q}}] \\[8pt]
&
\disp+ 3\|\mathcal{A}^*\hat{u}^k + \mathcal{B}^*\hat{v}^k - c\|^2\\[8pt]
\leq & 6\Phi_1(\bar{u}, \bar{v}, \bar{x})  + 3(2 + 8(\tau^{-1}-1)^2)\Phi_1(\bar{u},\bar{v}, \bar{x})+ 3C_4 + 3C_3,
\end{array}
\end{equation*}
where
$C_4  =   \max((\tau\sigma)^{-1}, \|\frac{1}{2}\mathcal{Q}+ \text{Diag}\,(\mathcal{D}_1 + \mathcal{S}+2\sigma \mathcal{A}\mathcal{A}^*, \mathcal{Q}_{22} + \mathcal{D}_2 + \mathcal{T})\|)$.
Then we can get that for any $(u,v,x)\in {B}_k$,  $\Phi_1(u,v,x) $ is bounded by a positive constant $$C_5 = (12 + 24(\tau^{-1}-1)^2)\Phi_1(\bar{u},\bar{v}, \bar{x})+ 3C_4 + 3C_3.$$ This, together with  (\ref{vi_complexity}), implies that
\begin{equation*}\begin{array}{ll}
&(p(\hat{u}^k) + q(\hat{v}^k) ) -(p(u) + q(v) )+\langle \hat{w}^{k}-w,\nabla\phi(w)\rangle+ \langle   \hat{u}^{k} - u, \mathcal{A}x\rangle  + \langle  \hat{v}^{k} - v , \mathcal{B}x\rangle\\[8pt]
&
 - \langle  \hat{x}^{k} - x , \mathcal{A}^*u + \mathcal{B}^*v - c\rangle
\leq  C_5/2k.
\end{array}
\end{equation*}
Under the conditions for case (b),
by a very similar approach we can get an upper bound $\overline{C}_5$ for $(\Psi_1(u,v,x) +\Xi_1)$
such that
\begin{equation*}\begin{array}{ll}
(p(\hat{u}^k) + q(\hat{v}^k)) -  (p(u) + q(v))  +\langle \hat{w}^{k}-w,\nabla\phi(w)\rangle+ \langle   \hat{u}^{k} - u, \mathcal{A}x\rangle  + \langle  \hat{v}^{k} - v , \mathcal{B}x\rangle\\[8pt]

- \langle  \hat{x}^k -x, \mathcal{A}^*u + \mathcal{B}^*v - c\rangle

 \leq (\Psi_1(u,v,x)+\Xi_1)/2k\leq \overline{C}_5/2k.
\end{array}
\end{equation*}
The above arguments show that for both cases,  the property (\ref{complexity_decrease_1}) holds.

\noindent
{\bf (iii)} For the complexity of primal objective functions, first,  we know from (\ref{kkt}) that
\begin{equation*}\begin{array}{ll}
p(u)\geq p(\bar{u}) + \langle -\mathcal{A}\bar{x} - \nabla_u\phi(\bar{w}), u-\bar{u}\rangle \quad \forall u\in\mathcal{U},\\[8pt]
q(v)\geq q(\bar{v}) + \langle -\mathcal{B}\bar{x} - \nabla_v\phi(\bar{w}), v - \bar{v}\rangle \quad\forall v\in\mathcal{V}.
\end{array}
\end{equation*}
Therefore, summing them up and by noting  $\mathcal{A}^*\bar{u} + \mathcal{B}^*\bar{v} = c$ and the convexity of function $\phi$, we have that
\begin{equation*}\begin{array}{ll}
\theta(u,v) - \theta(\bar{u}, \bar{v}) & \geq -\langle \bar{x}, \mathcal{A}^*u + \mathcal{B}^*v - c\rangle
 + \phi(w) - \phi(\bar{w}) - \langle \nabla\phi(\bar{w}), w-\bar{w}\rangle\\[8pt]
 &\geq -\langle \bar{x}, \mathcal{A}^*u + \mathcal{B}^*v - c\rangle
  \quad\forall u\in\mathcal{U}, v\in\mathcal{V}.
 \end{array}
\end{equation*}
Thus, with $(u,v) = (\hat{u}^{k},\hat{v}^{k})$, it holds that
\begin{equation}\label{ergfun_1}\begin{array}{ll}
\theta(\hat{u}^{k},\hat{v}^{k}) - \theta(\bar{u}, \bar{v}) & \geq -\langle \bar{x}, \mathcal{A}^*\hat{u}^{k} + \mathcal{B}^*\hat{v}^{k} - c\rangle
\\[8pt]
&
\geq \disp-\frac{1}{2}(\frac{1}{k}\|\bar{x}\|^2 + k\|\mathcal{A}^*\hat{u}^{k} + \mathcal{B}^*\hat{v}^{k} - c\|^2)\\[8pt]
&
\geq \disp-\frac{1}{2}(\|\bar{x}\|^2 + C_3)/k,
\end{array}
\end{equation}
where $C_3$ is the same constant as in (\ref{erg_1}).

For the reverse part,
by (\ref{phi_convex}) and (\ref{phi_majorize}) we can obtain that for any $i\geq 1$,
\begin{equation*}\begin{array}{ll}
\phi(w^{i+1}) \leq \disp\phi(w^{i})  + \langle \nabla \phi(w^{i}),w^{i+1} - w^i\rangle  + \frac{1}{2}\|w^{i+1} - w^i\|^2_{\mathcal{Q} + \mathcal{H}},\\[8pt]

\phi(\bar{w}) \geq \disp\phi(w^{i})  + \langle \nabla \phi(w^{i}),\bar{w}-w^i\rangle  + \frac{1}{2}\|\bar{w}-w^i\|^2_{\mathcal{Q}},
\end{array}
\end{equation*}
which,  indicate, that
\begin{equation}\label{erg_10}
\phi(w^{i+1}) - \phi(\bar{w})\leq \langle \nabla\phi(w^i), w^{i+1} - \bar{w}\rangle + \frac{1}{2}\|w^{i+1} - w^i\|^2_{\mathcal{Q} + \mathcal{H}} - \frac{1}{2}\|w^i - \bar{w}\|^2_{\mathcal{Q}}.
\end{equation}
Thus, (\ref{pqineq}) and (\ref{erg_10}) imply that for $\tau\in(0,1]$ and any $i\geq 1$,
\begin{equation}\label{ergfun_3}
\begin{array}{ll}
&\theta(u^{i+1}, v^{i+1}) -\theta(\bar u,\bar v)\\[8pt]

\leq & \disp\frac{1}{2}\|w^{i+1} - w^i\|^2_{\mathcal{Q}+\mathcal{H}} - \frac{1}{2}\|w^i-\bar{w}\|^2_{\mathcal{Q}} +\langle\bar{w} - w^{i+1}, \mathcal{Q}(w^{i+1}- w^i)\rangle + \langle \bar{u} - u^{i+1}, \mathcal{A}\tilde{x}^{i+1}\rangle \\[8pt]
&
\disp + \langle \bar{v} - v^{i+1}, \mathcal{B}\tilde{x}^{i+1}\rangle  - \langle  v^{i+1} - v^i, \mathcal{Q}_{12}^*(\bar{u} - u^{i+1})\rangle
 +\sigma \langle \mathcal{A}^*( u^{i+1}-\bar{u}), \mathcal{B}^*(v^{i+1}-v^i)\rangle\\[8pt]
 &
 \disp+\langle \bar{u} - u^{i+1}, ( \mathcal{D}_1 + \mathcal{S})(u^{i+1} - u^i)\rangle  + \langle  \bar{v} -v^{i+1}, ( \mathcal{D}_2+\mathcal{T})(v^{i+1} - v^i)\rangle\\[8pt]
\leq &\disp\frac{1}{2}(\Lambda_{i} - \Lambda_{i+1})
- \frac{1}{2}(\|u^{i+1} - u^i\|^2_{\mathcal{S}} +\|v^{i+1} - v^i\|^2_{\mathcal{T}} + \sigma\|\mathcal{A}^*u^{i+1} + \mathcal{B}^*v^i - c\|^2 \\[8pt]
&
\disp+\sigma(1-\tau)\|\mathcal{A}^*u^{i+1} + \mathcal{B}^*v^{i+1} - c\|^2) \\[8pt]
\leq  & \disp\frac{1}{2}(\Lambda_{i} - \Lambda_{i+1}).
\end{array}
\end{equation}
 Therefore, summing up the above inequalities over $i = 1, \cdots k$ and by using the convexity of function $\theta$  we can obtain that
\begin{equation}\label{ergfun_2}
 \theta(\hat{u}_k,\hat{v}_k) -  \theta(\bar{u}, \bar{v})
\leq (\Lambda_{1}(u,v,x) - \Lambda_{k+1}(u,v,x))/2k
\leq  \Lambda_1/2k.
\end{equation}
The inequalities (\ref{ergfun_1}) and (\ref{ergfun_2}) indicate that (\ref{erg_obj}) holds for case (a).

Next, assume that the conditions for case (b) hold.  Similar to (\ref{ergfun_3}),  we have that
\begin{equation*}
\begin{array}{ll}
&\theta(u^{i+1}, v^{i+1}) -\theta(\bar u,\bar v)\\[8pt]
\leq &\disp\frac{1}{2}(\overline{\Lambda}_{i} - \overline{\Lambda}_{i+1})
- \frac{1}{2}(\|u^{i+1} - u^i\|^2_{\mathcal{S}-\eta\mathcal{D}_1} +\|v^{i+1} - v^i\|^2_{\mathcal{T} + \min(\tau, 1+\tau - \tau^2)\sigma\mathcal{B}\mathcal{B}^* - \eta\mathcal{D}_2} \\[8pt]
&
+\min(1,1 + \tau^{-1}-\tau)\sigma\|\mathcal{A}^*u^{i+1} + \mathcal{B}^*v^{i+1} - c\|^2).
\end{array}
\end{equation*}
By the assumptions that $\mathcal{S}-\eta\mathcal{D}_1\succeq 0$ and $\mathcal{T} - \eta\mathcal{D}_2\succeq 0$, we can obtain that
\begin{equation}\label{ergfun_5}
\theta(\hat{u}^{k},\hat{v}^{k}) - \theta(\bar{u}, \bar{v})\leq (\overline{\Lambda}_1 - \overline{\Lambda}_{k+1})/2k\leq \overline{\Lambda}_1/2k.
\end{equation}
Thus, by (\ref{ergfun_2}) and (\ref{ergfun_5}) we can obtain the inequality (\ref{erg_obj}).
\qed
\bigskip

Below we make a couple of remarks about the results in Theorem \ref{erg_main}.
\begin{remark}
The result in part (ii) can be regarded as an ergodic version to (\ref{nonerg_main}) on the KKT optimality condition, though less explicit. Note that if one takes the square root on (\ref{nonerg_main}), the right hand side will become $O(1/\sqrt{k})$.
The inequality (\ref{complexity_decrease_1})
indicates that the majorized ADMM requires no more than $O(1/\varepsilon)$ iterations to obtain an $\varepsilon$-approximation solution in the sense of (\ref{def_approxiamteVI}).
When the objective function is separable, 
this kind of results 
has been studied for the (proximal) classic ADMM with separable objective functions, for examples,~\cite{Monteiro2013,HeYuan2012} and in a recent work by Li et al.~\cite{LiSunToh_indefinite2014}
 for a  majorized ADMM
%
with indefinite proximal terms.
  
\end{remark}

\begin{remark}
The results in parts (i) and (iii), which are on the ergodic complexity of the primal feasibility and the objective function, respectively, are extended from the work of Davis and Yin~\cite{Yin2014} on the classic ADMM with separable objective functions. These results are more explicit than the one in part (ii). However, there is no  corresponding result available on the dual problem. Therefore, it will be very interesting to see if one can develop a more explicit ergodic complexity result containing all the three parts in the KKT condition. 
\end{remark}

\section{Conclusions}
In this paper, we establish the convergence properties for the majorized ADMM with a large step length to solve linearly constrained convex programming whose {objective function} includes a coupled smooth function. From Theorem \ref{global_convergence}, one can see the influence of the coupled objective on the convergence condition. For $\tau\in(0,\frac{1+\sqrt{5}}{2})$, a joint condition like (\ref{case1_assumption}) or (\ref{case2_assumption}) is  needed
to analyze the behaviour of the iteration sequence. One can further observe that the parameter $\eta$, which controls the off-diagonal term of the generalized Hessian,  also affects the choice of proximal operators $\mathcal{S}$ and $\mathcal{T}$. However, as is pointed out in Remark \ref{Q_plus_separable}, when the {coupled function is convex quadratic}, $\eta = 0$ and the corresponding influence would disappear. Although, in this paper we focus on the 2-block case, it is not hard to see that, with the help of the Schur complement technique introduced in~\cite{LiSunToh_indefinite2014},
 one can apply our  majorized ADMM to solve  large scale convex optimization problems with many smooth blocks.

%
\bigskip

{\bf Acknowledgements.} The authors would like to thank Dr.~Caihua Chen at the Nanjing University for discussions on the iteration complexity described in the paper.

\end{document}